\documentclass{amsart}

\usepackage{amsthm}
\usepackage{amsfonts}

\newcommand{\C}{{\mathbb{C}}}
\newcommand{\Pp}{{\mathbb{P}}}

\newcommand{\Cy}[1]{X_{#1}}
\newcommand{\Inv}[1]{{\mathcal{I}}_{#1}}

\newtheorem{theorem}{Theorem}[section]
\newtheorem{corollary}[theorem]{Corollary}
\newtheorem{lemma}[theorem]{Lemma}
\newtheorem{proposition}[theorem]{Proposition}

\newtheorem{definition}[theorem]{Definition}

\title{Permutation representations on Schubert varieties}
\author{Julianna S. Tymoczko}
\address{Department of Mathematics, University of Iowa, Iowa City, IA 52242}
\email{tymoczko@math.uiowa.edu}
\thanks{MSC primary 20C30, secondary 05E10, 14M15, 55N91, 20F55. \\ \indent The author was partially supported by NSF grant 0402874.}

\begin{document}

\begin{abstract} This paper defines and studies permutation
representations on the 
equivariant cohomology of Schubert varieties, as representations both over
$\C$ and over $\C[t_1, t_2,\ldots,t_n]$.  We show these group actions are the same
as an action of simple transpositions studied geometrically by M.\ Brion, and give topological meaning to the divided difference operators of Berstein-Gelfand-Gelfand, Demazure, 
Kostant-Kumar, and others.  We analyze these representations using the 
 combinatorial approach to equivariant cohomology introduced by 
 Goresky-Kottwitz-MacPherson.  We find that each permutation
representation on equivariant cohomology produces a representation
on ordinary cohomology that is trivial, though the equivariant representation
 is not.  
\end{abstract}
\maketitle
\section{Introduction}

Geometric representation theory is an important approach to understanding permutation representations.  It builds algebraic varieties whose cohomology carries group actions and uses the geometric structure to analyze the representations. However, constructions of these representations to date are either not elementary (e.g. \cite{Sp}, \cite{BoM}, or \cite{L}) or not explicit (e.g. \cite{S} or \cite{P}).  We rectify this situation.  This paper constructs permutation representations on cohomology and equivariant cohomology of Schubert varieties in a simple yet concrete fashion.  Divided difference operators akin to those of Bernstein-Gelfand-Gelfand \cite{BGG} and Demazure \cite{D} result naturally from this geometric representation.

Denote the flag variety by the quotient $GL_n(\C)/B$, where $B$ is the group of invertible upper-triangular matrices.  Let $[g]$ be the flag corresponding to the matrix $g$.  For each permutation matrix $w$, the Schubert variety $\Cy{w}$ is the closure of the flags $[Bw]$ in $GL_n(\C)/B$. Schubert varieties are studied because they form a natural basis for the (equivariant) cohomology of the flag variety, and have deep combinatorial connections. The diagonal matrices in $GL_n(\C)$ act on each Schubert variety: if $t$ is diagonal and $[g]$ is in $\Cy{w}$ then $t \cdot [g] = [tg]$.  We study $H^*_T(\Cy{w})$ for this torus.  

This paper considers the action of the permutation group $S_n$ on equivariant and ordinary cohomology induced from the geometric action of $u \in S_n$ on $[g] \in GL_n(\C)/B$ given by $u \cdot [g] = [u^{-1}g]$.  Corollary \ref{W-action} provides a simple formula for the action of an {\em arbitrary} element of $S_n$ on each equivariant class in $H^*_T(GL_n(\C)/B)$.  We then obtain an explicit formula for the action of a simple transposition on the basis of equivariant Schubert classes, which is the core of our work.  If $s_{i,i+1} \in S_n$ is the simple transposition $(i,i+1)$ and $[\Omega_v]_{\Cy{w}}$ is a Schubert class in $\Cy{w}$ then
\begin{equation} \label{core step} s_{i,i+1}[\Omega_v]_{\Cy{w}} = \left\{ \begin{array}{ll}
[\Omega_v]_{\Cy{w}} & \textup{ if } s_{i,i+1}v > v \textup{ and} \\
{[\Omega_v]_{\Cy{w}}}+(t_{i+1}-t_i)[\Omega_{s_{i,i+1}v}]_{\Cy{w}} & \textup{ if } s_{i,i+1}v < v. \end{array} \right.\end{equation}
This formula was proven by M.\ Brion in \cite[Proposition 6.2]{B} using geometric methods that can only construct the action of simple transpositions on $H^*_T(\Cy{w})$. We show the action in \cite{B} derives from the global geometric action of $S_n$ on $GL_n(\C)/B$.

Equivariant cohomology surjects onto ordinary cohomology by the map that sends each $t_i \mapsto 0$.  Together with Equation \eqref{core step}, this is part of the main theorem:
\begin{theorem}
Fix a permutation $w$.  Denote the trivial representation in degree $d$ by $1^d$.  The action of $u \in S_n$ on $[g] \in GL_n(\C)/B$ given by $u \cdot [g] = [u^{-1}g]$ induces an action of $S_n$ on both $H^*_T(\Cy{w})$ and $H^*(\Cy{w})$ with the following properties:
\begin{enumerate}
\item $H^*_T(\Cy{w})$ is isomorphic to $\bigoplus_{[v] \in [S_n] \cap \Cy{w}}
1^{\ell(v)} \otimes \C[t_1, \ldots, t_n]$ as a graded $\C[S_n]$-module, where $S_n$ acts on 
$\C[t_1, \ldots, t_n]$ by permuting the variables;
\item $H^*_T(\Cy{w})$ is isomorphic to $\bigoplus_{[v] \in [S_n] \cap \Cy{w}}
1^{\ell(v)}$ as a graded twisted module over $\C[t_1,\ldots,t_n][S_n]$; and
\item $H^*(\Cy{w})$ is isomorphic to $\bigoplus_{[v] \in [S_n] \cap \Cy{w}}
1^{\ell(v)}$ as a graded $\C[S_n]$-module.
\end{enumerate}
\end{theorem}

Since the group $GL_n(\C)$ is connected, the endomorphism on $GL_n(\C)/B$ given by $u \cdot [g] = [u^{-1}g]$ is homotopic to the identity for each $u \in GL_n(\C)$, including $u \in S_n$.  Thus the action induced by $S_n$ on the ordinary cohomology $H^*(GL_n(\C)/B)$ is trivial.  However, the map $u \cdot [g] = [u^{-1}g]$ is not $T$-equivariant, so this fails for equivariant cohomology (indeed, Equation \eqref{core step} shows $S_n$ does not act trivially on equivariant Schubert classes).  Moreover, the action $u \cdot [g] = [u^{-1}g]$ is not generally well-defined on the Schubert variety $X_w$.  In fact $S_n$ does act on the cohomology of Schubert varieties and this action is trivial because it is a quotient of a direct sum of trivial representations.  An open question is to interpret this action geometrically.

Rewriting Equation \eqref{core step} gives the divided difference operator $[X] \mapsto \frac{[X]-s_{i,i+1}[X]}{t_i-t_{i+1}}$.  First defined by Bernstein-Gelfand-Gelfand \cite{BGG} and Demazure \cite{D}, divided difference operators have been widely used to analyze the algebraic structure of $H^*(GL_n(\C)/B)$
(e.g. \cite{KK1}, \cite{KK2}) and to find Schubert polynomials, namely ``nice" representatives for Schubert classes (e.g. \cite{LS}, \cite{BJS}, \cite{FK}); \cite{F2} has a survey.
This paper is atypical in treating $H^*_T(GL_n(\C)/B)$ as a {\em subring} of a product of $n!$ polynomial rings rather than as a quotient ring (see also \cite{KK1}, \cite{KK2}, and \cite{B}).  The group $S_n$ acts naturally on $H^*_T(GL_n(\C)/B)$ by both left and right multiplication on the index set, though this distinction is often elided in the literature.  The two actions give divided difference operators in very different ways.  We consider the action of left multiplication and its left divided difference operator; the other case is in \cite{T2}.  The formula for the left divided difference operator is that of Bernstein-Gelfand-Gelfand/Demazure, but the morphism is not.  This case is also studied in \cite{B}, where it is misidentified as the Bernstein-Gelfand-Gelfand/Demazure operator.  In fact, the right divided difference operators---defined by Kostant-Kumar in \cite{KK1}, \cite{KK2}---were proven in \cite{A} to be the divided difference operators of Bernstein-Gelfand-Gelfand and Demazure.    If $H^*_T(GL_n(\C)/B)$ is presented as a quotient ring, the left divided difference operator is the divided difference operator ``in the $y$-variables" of double Schubert polynomials (see e.g. \cite{F1}). Unusually, our arguments rely primarily on elementary combinatorics.

Our approach to equivariant cohomology uses the method of M.\ Goresky, R.\ Kottwitz, and R.\ MacPherson \cite{GKM}, which translates topological data of cohomology into a purely combinatorial calculation and is sketched in Section \ref{GKM theory}. (\cite{KK1} and \cite{KK2} essentially develop GKM theory by hand for $GL_n(\C)/B$.) The GKM method applies to Schubert varieties, as described in \cite{C}.   The ring $H^*_T(\Cy{w})$ has a special basis of ``Schubert classes"  produced by a combinatorial algorithm modeled on work of Guillemin-Zara and Knutson-Tao  \cite{GZ1}, \cite{GZ2}, \cite{KT}.  These classes are also the Kostant-Kumar $\xi^v$-classes of \cite{KK1} and \cite{KK2}. Section \ref{Knutson-Tao} discusses these bases for more general algebraic varieties than in \cite{GZ1}, \cite{GZ2}. Section \ref{eq reps} contains the equivariant cohomology calculations.

Our statements and proofs are first given for flags over $GL_n(\C)$.  Nonetheless, these results hold for all Lie types.  Our presentation was chosen because the exposition is more concrete for $GL_n(\C)$, because some readers will primarily be interested in this case, and because the reader interested in other Lie types can usually extrapolate those results immediately from our description of the special case.The general statements and all proofs that do not just change notation are in Section \ref{general type}.  

The author gratefully thanks Charles Cadman, William Fulton, Mel Hochster, Robert MacPherson, John Stembridge, and the referee for helpful comments.  
\section{Background}
\subsection{Permutation statistics}

We fix notation and give background about flag varieties and 
permutations.

As before, $B$ denotes the group of invertible upper-triangular matrices.
The flag variety is denoted $G/B$, and its typical element is denoted $[g]$.  
We also consider $[g]$ to be the collection of nested
subspaces whose $i$-dimensional subspace is the span of the first $i$ columns of $g$,
for each $i$.

Fix the standard basis $e_1$, $e_2$, $\ldots$, $e_n$ in $\C^n$.
Each permutation matrix $w$ also gives an element $[w]$ of the
flag variety.  We use the same notation for the matrix $w$ and for the permutation
on $\{1,2,\ldots,n\}$ given by $we_i = e_{w(i)}$.  We write  
$s_{jk}$ to denote the transposition that exchanges $j$ and $k$.

If $u$ is a matrix, its $(i,j)$ entry is denoted $u_{ij}$.  

\begin{definition}
For each permutation $w$, define the subgroup $U_w$ of $GL_n$ to be
\[\left\{u \in B: u_{ii}=1  \textup{ for each } i, \textup{ and if }i \neq j \textup{ then }
u_{ij}=0 \textup{ unless } w^{-1}(i) > w^{-1}(j)\right\}.\]
\end{definition}

The nonzero, nondiagonal entries in $U_w$ are important in what follows.

\begin{definition}
A pair $(i,j)$ is an inversion for $w$ if $i<j$ and $w^{-1}(i)>w^{-1}(j)$.  The
set $\Inv{w} = \{t_i-t_j: i<j, w^{-1}(i)>w^{-1}(j)\}$ is a collection of binomials
bijectively associated to the inversions for $w$.
\end{definition}

More generally, $\Inv{w}$ is the set of the positive roots that $w$ sends to negative roots.

\begin{definition}
For each $j \neq k$ and $c \in \C$, define the matrix $G_{jk}(c)$ by
\[\left( G_{jk}(c) \right)_{il} = \left\{ \begin{array}{rl} 1 & \textup{ if } i=l, \\
	c & \textup{ if } i=j \textup{ and }k=l, \textup{ and}\\
	0 & \textup{ otherwise.} \end{array} \right.\]
\end{definition}

The group $\{G_{jk}(c): c \in \C\}$
 is a subgroup of $U_w$ if $j<k$ and $w^{-1}(j)>w^{-1}(k)$.  Lie theoretically, 
the group $\{G_{jk}(c): c \in \C\}$ is a root subgroup.

If $w$ is a permutation, then $[Bw]$ is the Schubert cell corresponding to $w$
and $\Cy{w} = \overline{[Bw]}$ is the corresponding Schubert variety. 
For each $w$, the Schubert cell $[Bw]$ is homeomorphic to affine space and is parametrized
by the subgroup $U_w$. 
The flag variety is the disjoint union $\bigcup_{w \in S_n} [Bw]$.

\begin{proposition} \label{matrix rules}
Each flag in the Schubert cell $[Bw]$ can be written uniquely as a matrix of the form $U_ww$.
Each matrix $g \in U_ww$ can be written uniquely as $g= w+u$, where $u$ is 
zero except in entries $u_{ij}$ that are {\bf both} to the left {\bf and} above a nonzero
entry in $w$.
\end{proposition}

If $w$ and $v$ are permutations, we say that $w \geq v$ if $\overline{[Bw]} \supseteq
[Bv]$.  The {\em length} of the permutation $w$ is $\ell(w)=
\dim \overline{[Bw]}$.  We have $\ell(w)=| \Inv{w}|$ by construction.

The combinatorial argument in the next proof is similar to
\cite[Section 3]{KM}.  Similar results appear in the literature (especially \cite[Lemma 2.4]{BGG});
 since we were unable to find 
this formulation, we include it here. Section \ref{general type} contains the general proof.

\begin{lemma} \label{perm lemma}
 If $j<k$ and $w$ is a permutation with $\ell(s_{jk}w) = \ell(w)+1$ then 
	\begin{enumerate}
	\item \label{free} the inversions $\Inv{s_{jk}w} \cong 
		 \{t_j-t_k\} \cup \Inv{w} \textup{ mod}(t_j-t_k)$ with multiplicity; and 
	\item \label{transposing reflection} 
		if $s_{i,i+1}$ is a simple transposition with
		$s_{i,i+1}w > w$ then $s_{i,i+1}s_{jk}w > s_{jk}w$.
	\end{enumerate}
\end{lemma}

\begin{proof}
We use Proposition \ref{matrix rules} repeatedly.  Figure \ref{transposing flag} 
\begin{figure}[ht]
\[\left( \begin{array}{cc|c|ccc|c|cc}
	 \multicolumn{2}{c|}{\hspace{2em}} & A & \multicolumn{3}{|c|}{} & B & 
		\multicolumn{2}{|c}{\hspace{2em}} \\
	 \multicolumn{2}{c|}{\hspace{2em}} &  & \multicolumn{3}{|c|}{} &  & 
		\multicolumn{2}{|c}{\hspace{2em}} \\
	\cline{1-9} \multicolumn{2}{c|}{C} & 1 & \multicolumn{3}{|c|}{0 \cdots 0} & 0 & 
		\multicolumn{2}{|c}{0 \cdots 0} \\
	\cline{1-9} \multicolumn{2}{c|}{} & 0 & \multicolumn{3}{|c|}{} &  & \multicolumn{2}{|c}{} \\
	 \multicolumn{2}{c|}{} & \vdots & \multicolumn{3}{|c|}{} & D & \multicolumn{2}{|c}{} \\
	 \multicolumn{2}{c|}{} & 0 & \multicolumn{3}{|c|}{} &  & \multicolumn{2}{|c}{} \\
	\cline{1-9} \multicolumn{2}{c|}{E} & 0 & \multicolumn{3}{|c|}{F} & 1 & 
	\multicolumn{2}{|c}{0 \cdots 0} \\
	\cline{1-9} \multicolumn{2}{c|}{} & 0 & \multicolumn{3}{|c|}{} & 0 & \multicolumn{2}{|c}{} \\
	 \multicolumn{2}{c|}{} & \vdots & \multicolumn{3}{|c|}{} & \vdots & \multicolumn{2}{|c}{}  \\
	 \multicolumn{2}{c|}{} & 0 & \multicolumn{3}{|c|}{} & 0 & \multicolumn{2}{|c}{} 
	 	\end{array}\right)  \hspace{.5in} 
			\left( \begin{array}{cc|c|ccc|c|cc}
	 \multicolumn{2}{c|}{\hspace{2em}} & A & \multicolumn{3}{|c|}{} & B & 
		\multicolumn{2}{|c}{\hspace{2em}} \\
	 \multicolumn{2}{c|}{\hspace{2em}} &  & \multicolumn{3}{|c|}{} &  & 
		\multicolumn{2}{|c}{\hspace{2em}} \\
	\cline{1-9} \multicolumn{2}{c|}{C} & {\bf a} & \multicolumn{3}{|c|}{F'} & 1 & 
		\multicolumn{2}{|c}{0 \cdots 0} \\
	\cline{1-9} \multicolumn{2}{c|}{} &  & \multicolumn{3}{|c|}{} & 0 & \multicolumn{2}{|c}{} \\
	 \multicolumn{2}{c|}{} & D' & \multicolumn{3}{|c|}{} & \vdots & \multicolumn{2}{|c}{} \\
	 \multicolumn{2}{c|}{} &  & \multicolumn{3}{|c|}{} & 0 & \multicolumn{2}{|c}{} \\
	\cline{1-9} \multicolumn{2}{c|}{E} & 1 & \multicolumn{3}{|c|}{0 \cdots 0} & 0 & 
	\multicolumn{2}{|c}{0 \cdots 0} \\
	\cline{1-9} \multicolumn{2}{c|}{} & 0 & \multicolumn{3}{|c|}{} & 0 & \multicolumn{2}{|c}{} \\
	 \multicolumn{2}{c|}{} & \vdots & \multicolumn{3}{|c|}{} & \vdots & \multicolumn{2}{|c}{}  \\
	 \multicolumn{2}{c|}{} & 0 & \multicolumn{3}{|c|}{} & 0 & \multicolumn{2}{|c}{} 
	 	\end{array}\right)  \] 
		
		\hspace{.05in} $U_ww$ \hspace{2.1in} $U_{s_{jk}w}s_{jk}w$
\caption{Transposing $U_ww$} \label{transposing flag}
\end{figure}
is a schematic for  the matrices $U_{s_{jk}w}s_{jk}w$ and
$U_{w}w$.  We first confirm that the matrices are labeled correctly.
The matrices
differ only in the columns and rows
indicated in Figure \ref{transposing flag}.  
Regions $A$, $B$, $C$, and $E$ have the same free entries
in $U_ww$ as in $U_{s_{jk}w}s_{jk}w$.
The entry labeled ${\bf a}$ is free in one set but not in the other.  
Region $D'$ has all of the nonzero entries from $D$, plus perhaps additional
nonzero entries; similarly for $F'$ and $F$.  We conclude that the matrices
 on the right are translated
by a permutation of greater length than those on the left, which means they
must be $U_{s_{jk}w}s_{jk}w$.  


Now compare $U_w$ to $U_{s_{jk}w}$ in Figure \ref{transposing subgroups}.  
\begin{figure}[h]
\[\left( \begin{array}{cc|c|ccc|c|cc}
	 \multicolumn{2}{c|}{\hspace{2em}} & A & \multicolumn{3}{|c|}{} & B & 
		\multicolumn{2}{|c}{\hspace{2em}} \\
	 \multicolumn{2}{c|}{\hspace{2em}} &  & \multicolumn{3}{|c|}{} &  & 
		\multicolumn{2}{|c}{\hspace{2em}} \\
	\cline{1-9} \multicolumn{2}{c|}{0 \cdots 0} & 1 & \multicolumn{3}{|c|}{C} & 0 & 
		\multicolumn{2}{|c}{C} \\
	\cline{1-9} \multicolumn{2}{c|}{} & 0 & \multicolumn{3}{|c|}{} &  & \multicolumn{2}{|c}{} \\
	 \multicolumn{2}{c|}{} & \vdots & \multicolumn{3}{|c|}{} & D & \multicolumn{2}{|c}{} \\
	 \multicolumn{2}{c|}{} & 0 & \multicolumn{3}{|c|}{} &  & \multicolumn{2}{|c}{} \\
	\cline{1-9} \multicolumn{2}{c|}{0 \cdots 0} & 0 & \multicolumn{3}{|c|}{0 \cdots 0} & 1 & 
	\multicolumn{2}{|c}{E, F} \\
	\cline{1-9} \multicolumn{2}{c|}{} & 0 & \multicolumn{3}{|c|}{} & 0 & \multicolumn{2}{|c}{} \\
	 \multicolumn{2}{c|}{} & \vdots & \multicolumn{3}{|c|}{} & \vdots & \multicolumn{2}{|c}{}  \\
	 \multicolumn{2}{c|}{} & 0 & \multicolumn{3}{|c|}{} & 0 & \multicolumn{2}{|c}{} 
	 	\end{array}\right)  \hspace{.5in} 
			\left( \begin{array}{cc|c|ccc|c|cc}
	 \multicolumn{2}{c|}{\hspace{2em}} & B & \multicolumn{3}{|c|}{} & A & 
		\multicolumn{2}{|c}{\hspace{2em}} \\
	 \multicolumn{2}{c|}{\hspace{2em}} &  & \multicolumn{3}{|c|}{} &  & 
		\multicolumn{2}{|c}{\hspace{2em}} \\
	\cline{1-9} \multicolumn{2}{c|}{0 \cdots 0} & 1 & \multicolumn{3}{|c|}{C,F'} & {\bf a} & 
		\multicolumn{2}{|c}{C,F'} \\
	\cline{1-9} \multicolumn{2}{c|}{} & 0 & \multicolumn{3}{|c|}{} &  & \multicolumn{2}{|c}{} \\
	 \multicolumn{2}{c|}{} & \vdots & \multicolumn{3}{|c|}{} & D' & \multicolumn{2}{|c}{} \\
	 \multicolumn{2}{c|}{} & 0 & \multicolumn{3}{|c|}{} &  & \multicolumn{2}{|c}{} \\
	\cline{1-9} \multicolumn{2}{c|}{0 \cdots 0} & 0 & \multicolumn{3}{|c|}{0 \cdots 0} & 1 & 
	\multicolumn{2}{|c}{E} \\
	\cline{1-9} \multicolumn{2}{c|}{} & 0 & \multicolumn{3}{|c|}{} & 0 & \multicolumn{2}{|c}{} \\
	 \multicolumn{2}{c|}{} & \vdots & \multicolumn{3}{|c|}{} & \vdots & \multicolumn{2}{|c}{}  \\
	 \multicolumn{2}{c|}{} & 0 & \multicolumn{3}{|c|}{} & 0 & \multicolumn{2}{|c}{} 
	 	\end{array}\right)  \] 
		
		\hspace{.05in} $U_w$ \hspace{2.3in} $U_{s_{jk}w}$
\caption{The subgroups $U_w$ and $U_{s_{jk}w}$} \label{transposing subgroups}
\end{figure}
If $\ell(s_{jk}w) = \ell(w)+1$ then 
in fact $D$ and $D'$ have the same nonzero entries, as do $F$ and $F'$.
The entry marked ${\bf a}$ is free in 
$U_{s_{jk}w}$ but not in $U_w$, so $\Inv{s_{jk}w}$ contains
$t_j-t_k$ while $\Inv{w}$ does not.
Other than $t_j-t_k$, the multisets $\Inv{s_{jk}w} \textup{ mod}(t_j-t_k)$ and 
$\Inv{w} \textup{ mod}(t_j-t_k)$ are the same.  

For the last part, we also use Figure \ref{transposing subgroups}.  We are given that 
entry $(i,i+1)$ is zero in
$U_{w}$ and need to show that entry $(i,i+1)$ is zero in $U_{s_{jk}w}$.
If $(i,i+1)$ is in region $A$ in $U_w$ then entry $(i,k)$ in region $B$ of
$U_w$ is also zero, according to the description of the matrices $U_ww$ 
 in Proposition \ref{matrix rules}.  If $(i,i+1)$ is in region $C$ or $E$ in $U_w$ then 
it is unchanged in $U_{s_{jk}w}$.  This is true also if $(i,i+1)$ is in region $D$ of $U_w$,
since the free entries of $D$ and $D'$ are exactly the same. If $(i,i+1)$ is in region $F$, then entry
$(i,i+1)$ is zero in $U_{s_{jk}w}$ by construction.  If
$(i,i+1)$ is in region $F'$ then to be above the diagonal it must be between the $j^{th}$ and $k^{th}$
columns.  In that case, it is zero in $U_{s_{jk}w}$ since $F'$ and $F$ have the same free columns.
No other free entry in $U_w$ differs from that in $U_{s_{jk}w}$. 
\end{proof}

The next lemma specializes the previous argument.

\begin{lemma} \label{transposing inversions}
If $s_{i,i+1}w > w$, the inversions satisfy $\Inv{s_{i,i+1}w}=\{t_i-t_{i+1}\} \cup s_{i,i+1}\Inv{w}$,
	where $s_{i,i+1}$ acts on $\C[t_1,\ldots,t_n]$ by  
	$s_{i,i+1}(t_l)=t_{s_{i,i+1}(l)}$ for each $l$.
\end{lemma}

\begin{proof}
Consider Lemma \ref{perm lemma}.\ref{free} in the case when 
$s_{jk}=s_{i,i+1}$. The free entries in regions $C$ and $E$ of Figure
\ref{transposing flag} are exactly the same, and 
the regions marked $D$ and $D'$ do not exist. 
Figure \ref{transposing subgroups} now looks like:
\[ \left( \begin{array}{c|c|c|c}
& A & B &  \\
\cline{1-4} 0 \cdots 0 & 1 & 0 & C \\
\cline{1-4} 0 \cdots 0 & 0 & 1 & C, F \\
\cline{1-4} & \vdots & \vdots &
\end{array} \right)
\hspace{.5in}
 \left( \begin{array}{c|c|c|c}
 & B & A &  \\
\cline{1-4} 0 \cdots 0 & 1 & {\bf a} & C,F \\
\cline{1-4} 0 \cdots 0 & 0 & 1 & C \\
\cline{1-4} & \vdots & \vdots &
\end{array} \right) \]
The inversions in columns $i$ and $i+1$ and rows $i$ and $i+1$ are
exchanged.
\end{proof}

\subsection{GKM theory} \label{GKM theory}

The GKM method reduces the task
of identifying equivariant cohomology of a suitable space ${\mathcal X}$ to an algebraic
computation on a combinatorial graph associated to ${\mathcal X}$.  We 
sketch this method; the reader interested in more is encouraged
to see the original paper \cite{GKM}, the survey \cite{T1}, or the very nice
presentation in the introduction of \cite{KT}.

Let ${\mathcal X}$ be a complex projective algebraic variety with a linear algebraic action of 
a torus $T=\C^* \times \cdots \times \C^*$ that satisfies the following conditions: 
the torus has finitely 
many fixed points 
as well as finitely many one-dimensional orbits in ${\mathcal X}$; and
${\mathcal X}$ is {\em equivariantly formal}, a technical property that holds if, 
for instance, ${\mathcal X}$ has no odd-dimensional (ordinary) cohomology.  
These are the GKM conditions.  Any ${\mathcal X}$ satisfying the
GKM conditions is called a GKM variety.

The group $G$ acts on the flag variety $G/B$ by $h \cdot [g] = [hg]$.
If $T$ is the torus consisting of all diagonal matrices in $G$, 
the $G$-action restricts to a $T$-action.
Many subvarieties of the flag variety do not carry this torus action.  
However, Schubert varieties do
and are equivariantly formal with respect to this $T$-action.  We give the following
result for completeness, though it is not new.

\begin{proposition}
Every Schubert variety is $T$-equivariant and is equivariantly
formal with respect to this action.
\end{proposition}

\begin{proof}
Each Schubert variety has a cell decomposition as a union of Schubert cells.
Each Schubert variety is $T$-closed since it is $B$-closed.
No Schubert variety has odd-dimensional cohomology
since each $[Bw]$ is a complex affine cell.
Consequently, every Schubert variety is equivariantly formal with respect to every possible
torus action by \cite[Theorem 14]{GKM}.
\end{proof}

When ${\mathcal X}$ satisfies the GKM conditions, the closure of each
one-dimensional $T$-orbit in ${\mathcal X}$ is homeomorphic to $\Pp^1$.  
Both the origin and the point at 
infinity of $\Pp^1$ correspond to $T$-fixed points in ${\mathcal X}$.   If $O$ is a one-orbit in 
${\mathcal X}$ and 
$N_O$ and $S_O$ are the $T$-fixed points in $\overline{O}$, then 
the torus acts on the tangent space ${\mathcal T}_{N_O}(\overline{O})$ with weight $\alpha$ if and
only if the torus acts on the tangent space ${\mathcal T}_{S_O}(\overline{O})$ with weight $-\alpha$ (see
\cite[7.1.1]{GKM}).

In this case, we associate to ${\mathcal X}$ a labeled directed
graph called the
moment graph of ${\mathcal X}$.  The vertices of the moment graph of ${\mathcal X}$ are the
$T$-fixed points in ${\mathcal X}$, denoted ${\mathcal X}^T$.
If $v$ and $w$ are two vertices, there is an edge between $v$ and $w$ 
exactly when there is a 
one-orbit in ${\mathcal X}$ whose closure contains both $v$ and $w$.  Given a directed
edge $v \rightarrow w$ with associated one-orbit $O$, the edge $v \rightarrow w$ is labeled with the
weight of the torus action on the tangent space ${\mathcal T}_{v}(\overline{O})$.  (The
moment graph of ${\mathcal X}$ is not canonically directed; see Section \ref{Knutson-Tao} for more.)

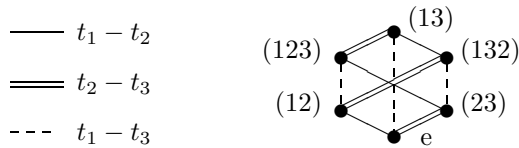
\begin{figure}[h]
\begin{picture}(350,60)(0,-25)
\put(175,25){\circle*{5}}
\put(155,15){\circle*{5}}
\put(195,15){\circle*{5}}
\put(155,-5){\circle*{5}}
\put(195,-5){\circle*{5}}
\put(175,-15){\circle*{5}}

\put(180,27){(13)}
\put(125,15){(123)}
\put(200,15){(132)}
\put(130,-5){(12)}
\put(200,-5){(23)}
\put(185,-18){e}

\put(175,25){\line(-2,-1){20}}
\put(175,27){\line(-2,-1){20}}
\put(175,-15){\line(2,1){20}}
\put(175,-13){\line(2,1){20}}
\put(195,15){\line(-2,-1){40}}
\put(195,17){\line(-2,-1){40}}

\put(175,25){\line(2,-1){20}}
\put(175,-15){\line(-2,1){20}}
\put(155,15){\line(2,-1){40}}

\put(195,15){\line(0,-1){5}}
\put(195,7){\line(0,-1){3}}
\put(195,1){\line(0,-1){3}}
\put(155,15){\line(0,-1){5}}
\put(155,7){\line(0,-1){3}}
\put(155,1){\line(0,-1){3}}
\put(175,25){\line(0,-1){5}}
\put(175,17){\line(0,-1){3}}
\put(175,11){\line(0,-1){3}}
\put(175,5){\line(0,-1){3}}
\put(175,-1){\line(0,-1){3}}
\put(175,-6){\line(0,-1){3}}
\put(175,-12){\line(0,-1){3}}

\put(30,25){\line(1,0){20}}
\put(55,22){$t_1-t_2$}
\put(30,-13){\line(1,0){3}}
\put(36,-13){\line(1,0){3}}
\put(42,-13){\line(1,0){3}}
\put(55,-16){$t_1-t_3$}
\put(30,6){\line(1,0){20}}
\put(30,4){\line(1,0){20}}
\put(55,3){$t_2-t_3$}

\end{picture}
\caption{The moment graph for $G/B$ when $n=3$} \label{moment graph eg}
\end{figure}

The flag variety is the example that forms the basis of the calculations in this paper.  
Part \ref{flag edges} of the next proposition describes the moment graph for the
flag variety $GL_n/B$, as stated (in more generality) in \cite[Theorem F]{C}.
Figure \ref{moment graph eg} gives the moment graph when $n=3$, drawn so that
each edge is directed from the higher endpoint to the lower endpoint.

Several corollaries follow immediately from the description of the moment graph.
Recall that if $V'$ is a subset of the vertices of the graph $(V,E)$, then the subgraph
{\em induced} by $V'$ is the maximal subgraph of $(V,E)$ with vertex set $V'$, namely
the graph $(V',E')$ where $E' = \{vw: v,w \in V' \textup{ and }vw \in E\}$.

\begin{proposition} \label{moment graph props}
In $GL_n(\C)/B$:
\begin{enumerate}
\item \label{flag edges} The vertices of the moment graph for $GL_n/B$ are
	$[w]$ for $w \in S_n$.  
	If $j<k$ then the flags $[G_{jk}(c)w]$ with $c \neq 0$	
	correspond to an edge directed out of $w$ if and only if $w^{-1}(j)>w^{-1}(k)$.  
	This edge is directed into $s_{jk}w$ and is labeled $t_j-t_k$.
\item There are $\ell(w)$ edges directed out of $w$, labeled exactly by the set $\Inv{w}$.  
\item If there is an edge from $s_{i,i+1}w$ to $w$ 
	and a transposition $s_{jk} \neq s_{i,i+1}$ with both 
		$\ell(s_{jk}w)=\ell(w)+1$ and an edge from $s_{jk}w$ to $w$, then
	there is also an edge from $s_{i,i+1}s_{jk}w$ to $s_{jk}w$.
\item \label{up neighbor} If $s_{i,i+1}w > w$ then the edges directed out of 
	$s_{i,i+1}w$ are labeled bijectively with $\{t_i-t_{i+1}\} \cup s_{i,i+1}\Inv{w}$,
	where $s_{i,i+1}$ acts on $\C[t_1,\ldots,t_n]$ by  the rule
	$s_{i,i+1}(t_l)=t_{s_{i,i+1}(l)}$ for each $l$.
\item \label{schubert moment graph}
	The moment graph for the Schubert variety $\Cy{w}$ is the subgraph of the
	moment graph for $GL_n/B$ induced by the permutation flags in $\Cy{w}$,
	namely $[v]$ where $v \in S_n$ satisfies $v \leq w$.
\end{enumerate}
\end{proposition}

\begin{proof}
Part \ref{flag edges} is \cite[Theorem F]{C}, with our convention for directing and labeling
the graph.  The number of edges directed out of $w$ is exactly
	$\dim \overline{[Bw]} = \ell(w)$, and is indexed by the inversions for $w$.  This proves
	Part 2.  Part 3 is Lemma \ref{perm lemma}.\ref{transposing reflection}.
Part 4 is Lemma \ref{transposing inversions}.  Part \ref{schubert moment graph} is 
in \cite[Theorem F5]{C}
and is a nice exercise for the reader.  \end{proof}

The main theorem that we use is:
\begin{proposition} {\bf (Goresky, Kottwitz, MacPherson)} 
Let ${\mathcal X}$ be a GKM variety and let $S$ be the polynomial ring $\C[t_1, \ldots, t_n]$.  
For each one-dimensional orbit $O$, the fixed points
in $\overline{O}$ are denoted $N_O$ and $S_O$, and the weight of the 
$T$-action on $O$ is $t_O$.
The equivariant cohomology $H^*_T({\mathcal X})$ is the subring of $S^{|{\mathcal X}^T|}$
given by:
\[H^*_T({\mathcal X}) = \left\{(p_w)_{w \in {\mathcal X}^T} \in S^{|{\mathcal X}^T|}:
 \textup{ for each one-orbit } O,  p_{N_O} - p_{S_O} \in \langle t_O \rangle \right\}.\]
\end{proposition}

An element $u \in S_n$ acts on the flag $[g] \in G/B$ by $u \cdot [g] = [u^{-1}g]$.  
This geometric action induces an action of $S_n$ on the equivariant cohomology of $G/B$ as
follows.  Recall
that $u \in S_n$ acts on the polynomial ring $\C[t_1, \ldots, t_n]$ by $u \cdot p(t_1,\ldots,t_n) =
p(t_{u(1)},\ldots,t_{u(n)})$.

\begin{corollary} \label{W-action}
The action of $S_n$ on $G/B$ induces a well-defined group action of $S_n$ 
on the equivariant cohomology $H^*_T(G/B)$, given by the rule that if $u \in S_n$
and $p = (p_v)_{v \in S_n} \in H^*_T(G/B)$ then the localization of $u \cdot p$ at each 
$v \in S_n$ is 
\[(u \cdot p)_v(t_1, \ldots,t_n) = u  \cdot p_{u^{-1}v}(t_1, \ldots, t_n) = p_{u^{-1}v}(t_{u(1)},\ldots, t_{u(n)}).\] 
\end{corollary}

\begin{proof}
The action of $u \in S_n$ on the variety $G/B$ gives a graph automorphism of the {\em undirected} 
moment graph  since $s_{u^{-1}(j),u^{-1}(k)}u^{-1}v=u^{-1}s_{jk}v$.  In particular, there is an edge between $v$ and $s_{jk}v$ if and only if there is an edge between $u^{-1}v$ and $u^{-1} s_{jk}v$. The edge between $u^{-1}v$ and $u^{-1}s_{jk}v$ is labeled $t_{u^{-1}(j)}-t_{u^{-1}(k)}$ or $t_{u^{-1}(k)}-t_{u^{-1}(j)}$.

This means the polynomial $(u \cdot p)_{v} - (u \cdot p)_{s_{jk}v}$ is in the ideal $\langle t_{j}-t_{k} \rangle$ if and only if $p_{u^{-1}v} - p_{u^{-1}s_{jk}v}$ is in $\langle t_{u^{-1}(j)}-t_{u^{-1}(k)} \rangle$.  So, the $S_n$-action is well-defined.

We now show this action is induced from the geometric action of $S_n$ on $G/B$.  We use the Borel construction of equivariant cohomology of a complex algebraic variety ${\mathcal X}$ with the action of a torus $T$.  If $ET$ is the classifying bundle of $T$, namely a contractible space on which $T$ acts freely, then the equivariant cohomology is defined to be $H^*_T({\mathcal X}) = H^*(ET \times^T {\mathcal X})$, where $ET \times^T {\mathcal X}$ denotes the quotient of the product $ET \times {\mathcal X}$ under the equivalence relation $(e,x) \sim (et, t^{-1}x)$ for all $t \in T$.  For the complex torus $T = \C^* \times \cdots \times \C^*$, one may use the classifying bundle $ET = (\C^{\infty} \backslash \{0\}) \times \cdots \times (\C^{\infty} \backslash \{0\})$.  The action of $u \in S_n$ on $(e,x) \in ET \times^T G/B$ is defined by $u \cdot (e,[g]) = (eu, [u^{-1}g])$, where $u$ acts on $e$ by permuting its coordinates.  Suppose $(et, [t^{-1}g]) \sim (e,[g])$.  There is a $t' \in T$ such that $(etw, [w^{-1}t^{-1}g]) = (ewt', [(t')^{-1}w^{-1}g])$ because the group $S_n$ normalizes the torus $T$.  Since $(ewt', [(t')^{-1}w^{-1}g]) \sim (ew, [w^{-1}g])$, this action of $S_n$ is well-defined on $ET \times^T G/B$.

The Borel construction of equivariant cohomology is related to the GKM construction by the inclusion $j: {\mathcal X}^T \hookrightarrow {\mathcal X}$, which for GKM spaces ${\mathcal X}$ induces an isomorphism onto its image $j^*: H^*_T({\mathcal X}) \hookrightarrow H^*_T({\mathcal X}^T)$.  The $S_n$-action on $ET \times^T G/B$ restricts to an $S_n$-action on $ET \times^T (G/B)^T$, which commutes with $j$ by definition.  Recall that 
\[H^*_T((G/B)^T) = \oplus_{w \in S_n} H^*_T(pt) = \oplus_{w \in S_n} \C[t_1, \ldots, t_n].\]  
For each $u \in S_n$, the map $u^*: H^*_T((G/B)^T) \rightarrow H^*_T((G/B)^T)$ permutes the fixed points by $u^{-1}$ and permutes the coordinates of the torus by $u$, which is precisely the action defined in this Corollary.  Since $j^* \circ u^* = u^* \circ j^*$, this is the action induced by the geometric action of $S_n$ on $G/B$.
\end{proof}

\begin{figure}[h]
\begin{picture}(350,60)(0,-25)
\put(125,25){\circle*{5}}
\put(105,15){\circle*{5}}
\put(145,15){\circle*{5}}
\put(105,-5){\circle*{5}}
\put(145,-5){\circle*{5}}
\put(125,-15){\circle*{5}}

\put(130,27){$t_1-t_3$}
\put(70,15){$t_1-t_2$}
\put(150,15){$t_1-t_3$}
\put(70,-5){$t_1-t_2$}
\put(150,-5){0}
\put(135,-18){0}

\put(125,25){\line(-2,-1){20}}
\put(125,27){\line(-2,-1){20}}
\put(125,-15){\line(2,1){20}}
\put(125,-13){\line(2,1){20}}
\put(145,15){\line(-2,-1){40}}
\put(145,17){\line(-2,-1){40}}

\put(125,25){\line(2,-1){20}}
\put(125,-15){\line(-2,1){20}}
\put(105,15){\line(2,-1){40}}

\put(145,15){\line(0,-1){5}}
\put(145,7){\line(0,-1){3}}
\put(145,1){\line(0,-1){3}}
\put(105,15){\line(0,-1){5}}
\put(105,7){\line(0,-1){3}}
\put(105,1){\line(0,-1){3}}
\put(125,25){\line(0,-1){5}}
\put(125,17){\line(0,-1){3}}
\put(125,11){\line(0,-1){3}}
\put(125,5){\line(0,-1){3}}
\put(125,-1){\line(0,-1){3}}
\put(125,-6){\line(0,-1){3}}
\put(125,-12){\line(0,-1){3}}

\put(275,25){\circle*{5}}
\put(255,15){\circle*{5}}
\put(295,15){\circle*{5}}
\put(255,-5){\circle*{5}}
\put(295,-5){\circle*{5}}
\put(275,-15){\circle*{5}}

\put(280,27){$t_2-t_3$}
\put(245,12){0}
\put(300,15){$t_2-t_3$}
\put(245,-8){0}
\put(300,-5){$t_2-t_1$}
\put(285,-18){$t_2-t_1$}

\put(275,25){\line(-2,-1){20}}
\put(275,27){\line(-2,-1){20}}
\put(275,-15){\line(2,1){20}}
\put(275,-13){\line(2,1){20}}
\put(295,15){\line(-2,-1){40}}
\put(295,17){\line(-2,-1){40}}

\put(275,25){\line(2,-1){20}}
\put(275,-15){\line(-2,1){20}}
\put(255,15){\line(2,-1){40}}

\put(295,15){\line(0,-1){5}}
\put(295,7){\line(0,-1){3}}
\put(295,1){\line(0,-1){3}}
\put(255,15){\line(0,-1){5}}
\put(255,7){\line(0,-1){3}}
\put(255,1){\line(0,-1){3}}
\put(275,25){\line(0,-1){5}}
\put(275,17){\line(0,-1){3}}
\put(275,11){\line(0,-1){3}}
\put(275,5){\line(0,-1){3}}
\put(275,-1){\line(0,-1){3}}
\put(275,-6){\line(0,-1){3}}
\put(275,-12){\line(0,-1){3}}

\put(25,2){$s_{1,2}$}
\put(45,2){$\cdot$}
\put(200,2){$=$}
\end{picture}
\caption{The permutation action on classes in $H^*_T(G/B)$} \label{examples of perm action}
\end{figure}
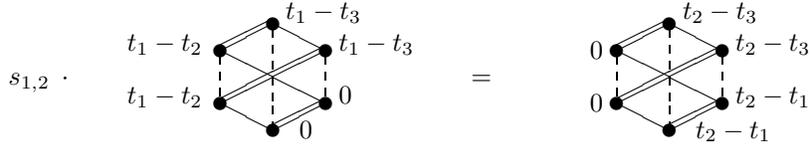

Figure \ref{examples of perm action} gives an example of this permutation
action on $H^*_T(G/B)$ when $n=3$.  
An equivariant 
class $p$ is denoted by a copy of the moment graph
in which each vertex $v$ is labeled with the polynomial $p_v$.  If we compute the 
action of $s_{2,3}$ instead of $s_{1,2}$ on the class in Figure \ref{examples of perm action},
 we obtain the same class we started with.

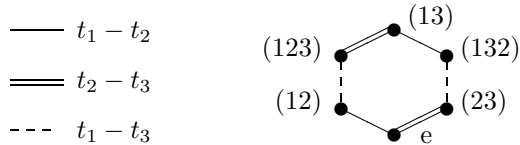
\begin{figure}[h]
\begin{picture}(350,60)(0,-25)
\put(175,25){\circle*{5}}
\put(155,15){\circle*{5}}
\put(195,15){\circle*{5}}
\put(155,-5){\circle*{5}}
\put(195,-5){\circle*{5}}
\put(175,-15){\circle*{5}}

\put(180,27){(13)}
\put(125,15){(123)}
\put(200,15){(132)}
\put(130,-5){(12)}
\put(200,-5){(23)}
\put(185,-18){e}

\put(175,25){\line(-2,-1){20}}
\put(175,27){\line(-2,-1){20}}
\put(175,-15){\line(2,1){20}}
\put(175,-13){\line(2,1){20}}

\put(175,25){\line(2,-1){20}}
\put(175,-15){\line(-2,1){20}}

\put(195,15){\line(0,-1){5}}
\put(195,7){\line(0,-1){3}}
\put(195,1){\line(0,-1){3}}
\put(155,15){\line(0,-1){5}}
\put(155,7){\line(0,-1){3}}
\put(155,1){\line(0,-1){3}}

\put(30,25){\line(1,0){20}}
\put(55,22){$t_1-t_2$}
\put(30,-13){\line(1,0){3}}
\put(36,-13){\line(1,0){3}}
\put(42,-13){\line(1,0){3}}
\put(55,-16){$t_1-t_3$}
\put(30,6){\line(1,0){20}}
\put(30,4){\line(1,0){20}}
\put(55,3){$t_2-t_3$}

\end{picture}
\caption{The toric variety for the decomposition into Weyl chambers in $GL_3$} \label{toric variety}
\end{figure}

In general, the moment graph of a subvariety of the flag variety 
is not an induced subgraph of the moment
graph for $G/B$.  For instance, Figure \ref{toric variety}
shows
the moment graph for the toric variety associated to the decomposition into Weyl chambers
in $GL_3$.  The (ordinary) cohomology of this toric variety carries a nontrivial $S_3$-action
that was studied in \cite{S} and \cite{P}.

\subsection{A combinatorial basis for $H^*_T(\mathcal{X})$}  \label{Knutson-Tao}

Once the moment graph is constructed, the GKM description of equivariant cohomology is 
purely combinatorial.  We now discuss a combinatorial construction of a
 basis for the equivariant cohomology ring as a module over $H^*_T(\text{pt})$,
which we call Knutson-Tao classes after \cite{KT}.  Our results are similar to those of 
Guillemin-Zara in \cite{GZ1} and \cite{GZ2} but use an underlying combinatorial model
that is simpler and less restrictive.  

\begin{definition}
A combinatorial moment graph is a finite graph whose edges are labeled by linear forms and are
directed so that there are no directed circuits.
\end{definition}
We will assume one additional condition on moment graphs in this section, which
holds for moment graphs of GKM varieties for geometric reasons:
\begin{itemize}
\item For each $v \in {\mathcal X}^T$, if 
$\beta_1, \ldots, \beta_k$ label the edges directed out of $v$ then 
the $\beta_i$ are pairwise linearly independent.
\end{itemize}
It is a small exercise to see that that is equivalent to the following condition:
\begin{itemize}
\item For each $v \in {\mathcal X}^T$, if 
$\beta_1, \ldots, \beta_k$ label the edges directed out of $v$ then 
$\langle \beta_1 \beta_2 \cdots \beta_k \rangle = \bigcap_{i=1}^k \langle \beta_i \rangle$.
\end{itemize}
If ${\mathcal X}$ is a GKM variety, then its moment graph can be directed without circuits
by choosing a suitably generic one-dimensional
subtorus $T' \subseteq T$ and directing the edges according to the
flow of $T'$ (see \cite[Section 5]{T1}).  Once its moment graph is directed acyclically, 
each edge $v \rightarrow w$ is labeled 
with the $T$-weight at $v$ on the corresponding one-dimensional orbit.  This gives
(many) combinatorial moment graphs associated to ${\mathcal X}$.

Combinatorial moment graphs form a larger class of graphs than studied by Guillemin-Zara, 
including for instance moment graphs of singular GKM varieties.   
We note that combinatorial moment graphs need not be regular and have no 
particular relationship between edges directed in and out of each vertex 
(axiom A3 of the Guillemin-Zara axial function,
 \cite[Definition 2.1.1]{GZ1} or \cite[Definition 2.1]{GZ2}).  

Using combinatorial moment graphs, 
we will show that Knutson-Tao classes are unique for a large class of varieties including
Schubert varieties.  For Grassmannians and flag varieties, the Knutsen-Tao classes are
Schubert classes.  We ask whether the results of Guillemin-Zara can be extended to show
that Knutson-Tao classes exist for combinatorial moment graphs.

Write $u \succeq_D u'$ if there is a directed path from $u$ to $u'$ in the combinatorial
moment graph, possibly of length zero.  For any  directed graph with no directed circuits, the
relation $\succeq_D$ is a partial order on the vertices.  
(This partial order coincides with the 
Bruhat order for the  moment
graph of $G/B$ in Proposition \ref{moment graph props}.)

We now define Knutson-Tao classes, which T.\ Braden and R.\ MacPherson also used
to construct equivariant intersection cohomology \cite{BrM} and which are homogeneous
versions of the generating classes in \cite[Definition 2.3]{GZ2}.

\begin{definition}
Let ${\mathcal X}$ be a GKM variety 
and let $v$ be a $T$-fixed point in ${\mathcal X}$.  A Knutson-Tao class for $v$ is
an equivariant class $(p^v_w)_{w \in {\mathcal X}^T} \in H^*_T({\mathcal X})$
for which:
\begin{enumerate}
\item $\displaystyle p^v_v = \prod_{\scriptsize \begin{array}{c} w \in {\mathcal X}^T \textup{ s.t. }\\ v \rightarrow w \textup{ is an edge}
\end{array}} \alpha_{vw}$, where $\alpha_{vw}$ is the label on the edge $v \rightarrow w$;
\item each nonzero $p^v_w$ is a homogeneous polynomial with $\deg(p^v_w)=\deg(p^v_v)$; and 
\item $p^v_w = 0$ for each $w$ with no directed path from $w$ to $v$, that is if $w \not \succ_D v$.
\end{enumerate}
\end{definition}

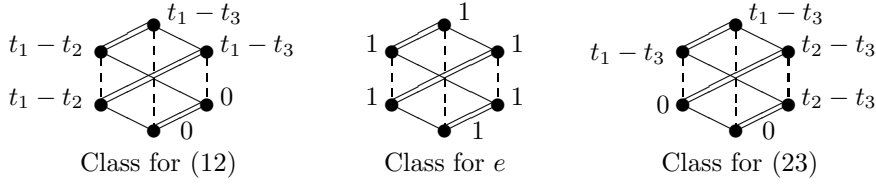
\begin{figure}[h]
\begin{picture}(350,65)(-5,-30)
\put(55,25){\circle*{5}}
\put(35,15){\circle*{5}}
\put(75,15){\circle*{5}}
\put(35,-5){\circle*{5}}
\put(75,-5){\circle*{5}}
\put(55,-15){\circle*{5}}

\put(60,27){$t_1-t_3$}
\put(0,15){$t_1-t_2$}
\put(80,15){$t_1-t_3$}
\put(0,-5){$t_1-t_2$}
\put(80,-5){0}
\put(65,-18){0}

\put(55,25){\line(-2,-1){20}}
\put(55,27){\line(-2,-1){20}}
\put(55,-15){\line(2,1){20}}
\put(55,-13){\line(2,1){20}}
\put(75,15){\line(-2,-1){40}}
\put(75,17){\line(-2,-1){40}}

\put(55,25){\line(2,-1){20}}
\put(55,-15){\line(-2,1){20}}
\put(35,15){\line(2,-1){40}}

\put(75,15){\line(0,-1){5}}
\put(75,7){\line(0,-1){3}}
\put(75,1){\line(0,-1){3}}
\put(35,15){\line(0,-1){5}}
\put(35,7){\line(0,-1){3}}
\put(35,1){\line(0,-1){3}}
\put(55,25){\line(0,-1){5}}
\put(55,17){\line(0,-1){3}}
\put(55,11){\line(0,-1){3}}
\put(55,5){\line(0,-1){3}}
\put(55,-1){\line(0,-1){3}}
\put(55,-6){\line(0,-1){3}}
\put(55,-12){\line(0,-1){3}}

\put(165,25){\circle*{5}}
\put(145,15){\circle*{5}}
\put(185,15){\circle*{5}}
\put(145,-5){\circle*{5}}
\put(185,-5){\circle*{5}}
\put(165,-15){\circle*{5}}

\put(170,27){1}
\put(135,15){1}
\put(190,15){1}
\put(135,-5){1}
\put(190,-5){1}
\put(175,-18){1}

\put(165,25){\line(-2,-1){20}}
\put(165,27){\line(-2,-1){20}}
\put(165,-15){\line(2,1){20}}
\put(165,-13){\line(2,1){20}}
\put(185,15){\line(-2,-1){40}}
\put(185,17){\line(-2,-1){40}}

\put(165,25){\line(2,-1){20}}
\put(165,-15){\line(-2,1){20}}
\put(145,15){\line(2,-1){40}}

\put(185,15){\line(0,-1){5}}
\put(185,7){\line(0,-1){3}}
\put(185,1){\line(0,-1){3}}
\put(145,15){\line(0,-1){5}}
\put(145,7){\line(0,-1){3}}
\put(145,1){\line(0,-1){3}}
\put(165,25){\line(0,-1){5}}
\put(165,17){\line(0,-1){3}}
\put(165,11){\line(0,-1){3}}
\put(165,5){\line(0,-1){3}}
\put(165,-1){\line(0,-1){3}}
\put(165,-6){\line(0,-1){3}}
\put(165,-12){\line(0,-1){3}}

\put(275,25){\circle*{5}}
\put(255,15){\circle*{5}}
\put(295,15){\circle*{5}}
\put(255,-5){\circle*{5}}
\put(295,-5){\circle*{5}}
\put(275,-15){\circle*{5}}

\put(280,27){$t_1-t_3$}
\put(220,12){$t_1-t_3$}
\put(300,15){$t_2-t_3$}
\put(245,-8){0}
\put(300,-5){$t_2-t_3$}
\put(285,-18){$0$}

\put(275,25){\line(-2,-1){20}}
\put(275,27){\line(-2,-1){20}}
\put(275,-15){\line(2,1){20}}
\put(275,-13){\line(2,1){20}}
\put(295,15){\line(-2,-1){40}}
\put(295,17){\line(-2,-1){40}}

\put(275,25){\line(2,-1){20}}
\put(275,-15){\line(-2,1){20}}
\put(255,15){\line(2,-1){40}}

\put(295,15){\line(0,-1){5}}
\put(295,7){\line(0,-1){3}}
\put(295,1){\line(0,-1){3}}
\put(255,15){\line(0,-1){5}}
\put(255,7){\line(0,-1){3}}
\put(255,1){\line(0,-1){3}}
\put(275,25){\line(0,-1){5}}
\put(275,17){\line(0,-1){3}}
\put(275,11){\line(0,-1){3}}
\put(275,5){\line(0,-1){3}}
\put(275,-1){\line(0,-1){3}}
\put(275,-6){\line(0,-1){3}}
\put(275,-12){\line(0,-1){3}}

\put(27,-30){Class for $(12)$}
\put(142,-30){Class for $e$}
\put(247,-30){Class for $(23)$}
\end{picture}
\caption{Three Knutson-Tao classes in $H^*_T(G/B)$} \label{egs of KT classes}
\end{figure}

Figure \ref{egs of KT classes} gives three examples of Knutson-Tao classes in $G/B$.  
The reader may note
that the Knutson-Tao class for the permutation $(12)$ appeared in Figure \ref{examples of perm action}.

When they exist, Knutson-Tao classes form a basis for equivariant cohomology. The following is similar to \cite[Theorem 2.4.4]{GZ1}.  The proof uses minimal elements in the moment graph; these exist because the graph is finite, directed, and acyclic.

\begin{proposition}
Suppose that for each $v \in {\mathcal X}^T$ there exists at least one Knutson-Tao class 
$p^v \in H^*_T({\mathcal X})$.  Then the classes $\{ p^v: v \in {\mathcal X}^T\}$ form a basis for $H^*_T({\mathcal X})$.
\end{proposition}

\begin{proof}
The proof will be by induction on the partial order defined by the moment graph.  For each $q \in H^*_T({\mathcal X})$ and each minimal element $v$ in the poset, we have $q_v = c_v p^v_v$ for a unique polynomial $c_v$.  (In fact $c_v = q_v$.)  Thus $(q - c_vp^v)_v = 0$.  

Let $S$ be any subset of ${\mathcal X}^T$ such that for each $v \in S$ and $u \prec_D v$ in the moment graph, the fixed point $u \in S$ as well.  The inductive hypothesis is that for each $q \in H^*_T({\mathcal X})$ there are unique coefficients $c_v$ such that $(q- \sum_{v \in S} c_vp^v)_u= 0$ for each $u \in S$.  For each minimal element $u' \in {\mathcal X}^T - S$, note that $(q- \sum_{v \in S} c_vp^v)_{u}= 0$ for all $u \prec_D u'$.  The GKM conditions imply there is a unique polynomial $c_{u'}$ with
\[(q- \sum_{v \in S} c_vp^v)_{u'} = c_{u'} p^{u'}_{u'}.\] 
So $(q- \sum_{v \in S} c_vp^v - c_{u'}p^{u'})_{u'} =0$.  Since $p^{u'}_{u} = 0$ for all $u \not \succ_D u'$, it is also true that $(q- \sum_{v \in S} c_vp^v - c_{u'}p^{u'})_{u} =(q- \sum_{v \in S} c_vp^v)_u = 0$ for all $u \in S$.  Thus the inductive hypothesis holds for the set $S \cup \{u'\}$.  There are a finite number of fixed points in ${\mathcal X}^T$ so there is a unique expression $q = \sum_{v \in {\mathcal X}^T} c_vp^v$ for each $q \in H^*_T({\mathcal X})$.
\end{proof}

Our next results hold for a family of varieties called 
{\em Palais-Smale} varieties.
Motivated by the comment in \cite[page 187]{K}, we generalize the definition of Palais-Smale
so that it applies to varieties that are not Hamiltonian, nor even smooth.

\begin{definition}
The GKM variety ${\mathcal X}$
 is Palais-Smale if its moment graph  can be directed 
so that there is an edge from $v$ to $u$ only if there are more edges directed out of $v$
 than out of $u$.
\end{definition}

For instance, each projective space $\C\Pp^n$ is Palais-Smale.  We will see that each Schubert
variety in the full flag variety is also Palais-Smale.  The toric variety whose moment graph is shown in 
Figure \ref{toric variety} is {\em not} Palais-Smale.

\begin{lemma}
Each Schubert variety is Palais-Smale, including $G/B$.
\end{lemma}

\begin{proof}
Our convention for the moment graph for $\Cy{w}$ is that there is a directed path
from $u$ to $v$ only if $u > v$ in the Bruhat order.  A directed circuit is a directed path from $u$ to $u$.
Since $u \not > u$, no such circuit exists in the moment graph for $\Cy{w}$.

There is a 
	directed edge from $s_{jk}v$ to $v$ only if $s_{jk}v>v$. 	
	By Proposition \ref{moment graph props}.\ref{flag edges}, there are $\ell(v)$ edges directed out of
	$v$ in $G/B$ and $\ell(s_{jk}v)$ edges out of $s_{jk}v$.
	This holds in $\Cy{w}$ by  Proposition \ref{moment graph props}.\ref{schubert moment graph}.
	If $s_{jk}v>v$ then by definition $\ell(s_{jk}v)>\ell(v)$.
\end{proof}

If ${\mathcal X}$ is Palais-Smale, then Knutson-Tao classes are unique.   The next
proposition is similar to \cite[Theorem 2.3]{GZ2}.  Like theirs,
it could be modified to prove uniqueness of a Knutson-Tao
class at one particular vertex.

\begin{lemma} \label{uniqueness of KT classes}
If ${\mathcal X}$ is Palais-Smale and $p^v=(p^v_w)_{w \in {\mathcal X}^T}$, $q^v=(q^v_w)_{w \in {\mathcal X}^T}$ are two Knutson-Tao
classes corresponding to $v$, then $p^v_w=q^v_w$ for each $w \in {\mathcal X}^T$.
\end{lemma}

\begin{proof}
The proof mimics \cite[Lemma 1]{KT}.  
Consider the class $(p^v-q^v) \in H^*_T({\mathcal X})$.  We know that $(p^v-q^v)_u=0$ if $u=v$ or
if there is no directed path from $u$ to $v$.  
Choose a minimal vertex $u_0$ satisfying $(p^v-q^v)_{u_0} \neq 0$.  This means both
that $u_0 \succ_D v$ and that $(p^v-q^v)_{u'}=0$
for all $u'$ with $u_0 \succ_D u'$. In particular, the polynomial $(p^v-q^v)_{u_0}$ is in the ideal
generated by the labels on the edges $u_0 \rightarrow u'$ by the GKM rules.  
The variety ${\mathcal X}$ is 
Palais-Smale so the number of edges $u_0 \rightarrow u'$ is greater than $\deg p^v$.
Since $p^v-q^v$ has degree $\deg p^v$, we conclude $(p^v-q^v)_{u_0} = 0$.
\end{proof}

It is not a priori clear that {\em any} Knutson-Tao classes exist.  We prove existence 
for $G/B$ in Section \ref{eq reps}.  Existence for various families of smooth GKM varieties
is proven in \cite{GZ2}.

\section{Permutation representations on the \\ (equivariant) cohomology of $\Cy{w}$} \label{eq reps}

In this section we study permutation representations on $H^*_T(\Cy{w})$ and
$H^*(\Cy{w})$.
In Section \ref{schubert classes}, we show how to construct a basis of 
Knutson-Tao classes for $H^*_T(\Cy{w})$ by restricting the Schubert basis for $H^*_T(G/B)$.  In 
Section \ref{flag repns}, we explicitly identify the $S_n$-action on $H^*_T(G/B)$ by computing how
each simple transposition acts on each Schubert class.  Finally, in
Section \ref{decomposing repns}, we use the formula of Section \ref{flag repns} and
a restriction map $\iota: H^*_T(G/B) \rightarrow H^*_T(\Cy{w})$ to analyze the 
restricted $S_n$-representation
on $H^*_T(\Cy{w})$ and $H^*(\Cy{w})$.

\subsection{Knutson-Tao classes in $\Cy{w}$ and $G/B$} \label{schubert classes}

\begin{lemma}
Define the map 
$\iota: H^*_T(G/B) \rightarrow H^*_T(\Cy{w})$ by restriction:
\[\iota\left((p_v)_{v \in S_n}\right) = \left( p_v \right)_{[v] \in \Cy{w}}.\]
The map $\iota$ is a well-defined ring and $\C[t_1,\ldots,t_n]$-module homomorphism.  
\end{lemma}

\begin{proof}
The moment graph of $\Cy{w}$ is an induced subgraph
of the moment graph of $G/B$, so for each $p \in H^*_T(G/B)$,
the GKM conditions hold on $\iota(p)$.  This means $\iota(p) \in H^*_T(\Cy{w})$.
 By construction, $\iota$ is both a ring and module homomorphism.
\end{proof}

We use the following property of equivariant cohomology \cite[Fact 2, page 10]{KT}: 
Let $X$ be a $T$-invariant oriented cycle in ${\mathcal X}$, a smooth compact complex algebraic variety.
Then $X$ determines a class $[X] \in H^*_T({\mathcal X})$ whose degree 
is the codimension of $X$ in ${\mathcal X}$.  If $u \in {\mathcal X}^T$ is not in $X$, then the localization of $[X]$ at $u$ is zero.  

For instance, let $[v] \in S_n$ be a permutation flag and let $B^-$ be the group of lower-triangular
invertible matrices.  The closure $\overline{[B^-vB]}$ in $G/B$ determines
a class $[\Omega_v] \in H^*_T(G/B)$ called the Schubert class of $v$ in $G/B$.

We give several properties of $[\Omega_v]$, which is also the class $\xi^v$ in \cite{KK1}.

\begin{lemma} \label{localizing schubert classes}
Write the class $[\Omega_v] \in H^*_T(G/B)$ as $[\Omega_v]=(p_u^v)_{u \in S_n}$.
\begin{enumerate}
\item \label{support} The $T$-fixed points in $\overline{[B^-v]}$ 
	are exactly those $u \in S_n$ with $u \geq v$.
\item \label{degree} If $v,u \in S_n$ and $u \geq v$ then the degree of $p_u^v$ is $\ell(v)$.
	If $u \not \geq v$ then $p_u^v=0$.
\item \label{localization of class} For each $v \in S_n$, the polynomial $p_v^v$ satisfies
	\[p_v^v= \prod_{t_i-t_j \in \Inv{v}} (t_i - t_j).\]
\item \label{Knutson-Tao class} The class $[\Omega_v]$ is the Knutson-Tao class corresponding
	to $v$ in $H^*_T(G/B)$.
\item \label{neighbors} Suppose $u>v$ has $\ell(u) = \ell(v)+1$ and that 
	the edge from $u$ to $v$ is labeled
	$t_{j} - t_{k}$.  Then the polynomial $p_u^v$ satisfies
	\[p_u^v= \prod_{\scriptsize \begin{array}{c}t_i-t_{i'} \in \Inv{u} \\ t_i - t_{i'} \neq t_j-t_k \end{array}} (t_i - t_{i'}).\]
\end{enumerate}
\end{lemma}

\begin{proof}
\begin{enumerate}
\item The closure relation $\overline{[B^-u]} \subseteq \overline{[B^-v]}$ is equivalent to
$u \geq v$.  
\item We use \cite[Fact 2, page 10]{KT}. 
The degree of $[\Omega_v]$ is the codimension of $\overline{[B^-v]}$ in $G/B$, 
namely $\ell(v)$.
If $u \not \geq v$ then $u \not \in \overline{[B^-v]}$ and so $p_u^v = 0$.
\item The localization $p_v^v$ is the product of the weights of the torus action 
on the normal space at $[v]$ to
$[\Omega_v]$ in $G/B$ (see \cite[Fact 3, p. 10]{KT}).
The normal space at $v$ to $\overline{[B^-v]}$ in $G/B$ is $[U_vv]$
by, for instance,  \cite[Theorem G.2]{C}.
\item The relation $u \not \geq v$ means there is no
path from $u$ to $v$ in the moment graph of $G/B$.  Parts \ref{degree} and
\ref{localization of class} show that $[\Omega_v]$ is a Knutson-Tao class.   
Knutson-Tao classes are unique by Lemma 
\ref{uniqueness of KT classes}, since $G/B$ is Palais-Smale.
\item  Let $q =\prod_{\scriptsize \begin{array}{c}t_i-t_{i'} \in \Inv{u} \\ t_i - t_{i'} \neq t_j-t_k \end{array}} (t_i - t_{i'})$.
 There are $\ell(v)+1$ edges directed out of $u$.  One 
edge is to $v$, so there are $\ell(v)$ edges from $u$ to vertices $u'$
with $p_{u'}^v = 0$.  The polynomial $p_u^v$ is in the ideal generated by the
labels of these $\ell(v)$ edges, so $p_u^v$ is a scalar multiple of $q$.  
Since $u=s_{jk}v$, we have 
$\Inv{v} \cong \{0\} \cup \Inv{s_{jk}v} \textup{ mod}(t_j-t_k)$ 
by Lemma \ref{perm lemma}.\ref{free}.  By Part \ref{localization of class},
the polynomial
$q - p_v^v$ is in $\langle t_j - t_k \rangle$. 
Since $q \not \in \langle t_j-t_k \rangle$, 
we conclude $p_u^v=q$.\end{enumerate}\end{proof}

The next lemma 
permits us to specialize calculations from $G/B$ to $\Cy{w}$.

\begin{lemma}
For each $[v] \in \Cy{w}$, restriction preserves Knutson-Tao classes:
\[ [\Omega_v]_{\Cy{w}} := \iota \left([\Omega_v]\right)  \textup{ is the unique Knutson-Tao class for $v$ in $\Cy{w}$.}\]
\end{lemma}

\begin{proof}
We
need to check that $\iota([\Omega_v])$ satisfies the Knutson-Tao conditions.

The map $\iota$ preserves degree, so $\deg [\Omega_v]_{\Cy{w}}$ is $\ell(v)$.
If $[u] \in \Cy{w}$ satisfies $u \not \geq v$ in $\Cy{w}$ then $u \not \geq v$ in 
$G/B$, so neither the moment graph for $\Cy{w}$ nor for $G/B$ has a path from $u$ to $v$.
For each such $u$, the localization of $[\Omega_v]_{\Cy{w}}$ to $u$ is zero
by definition of $\iota$.   
Lemma \ref{localizing schubert classes}.\ref{localization of class} proves 
the localization of $[\Omega_v]_{\Cy{w}}$ to $v$ is as desired.  The Knutson-Tao class is unique
because $\Cy{w}$ is Palais-Smale (Lemma \ref{uniqueness of KT classes}).
\end{proof}

When $\Cy{w}$ is singular, the Knutson-Tao class for $v$ need not be the geometric Schubert
class corrresponding to $v$.  The precise 
relationship between geometric Schubert classes 
and Knutson-Tao classes in $H^*_T(\Cy{w})$ is unknown.

The next corollary combines Lemma \ref{moment graph props}.\ref{up neighbor} 
with Parts \ref{localization of class}
and \ref{neighbors} of Lemma \ref{localizing schubert classes}.  It is the heart of our later calculations.

\begin{corollary} \label{neighboring polys}
Let $v, s_{i,i+1} \in S_n$ satisfy $s_{i,i+1}v>v$.   Write
$[\Omega_v]_{G/B}=(p_u^v)_{u \in S_n}$ and $[\Omega_{s_{i,i+1}v}]_{G/B}=
(p_u^{s_{i,i+1}v})_{u \in S_n}$ for the corresponding Schubert classes.  Then
\[ s_{i,i+1}p_{s_{i,i+1}v}^v = p_v^v = \frac{s_{i,i+1}p_{s_{i,i+1}v}^{s_{i,i+1}v}}{t_{i+1}-t_i}.\]
\end{corollary}

\subsection{The $S_n$-action on $H^*_T(G/B)$} \label{flag repns}

Recall that if $u \in S_n$ then $u$ acts on $p = (p_v(t_1, \ldots, t_n))_{v \in S_n} \in H^*_T(G/B)$
 by
\[(u \cdot p)_v = u \cdot p_{u^{-1}v}(t_1,\ldots,t_n) = p_{u^{-1}v}(t_{u(1)},\ldots,t_{u(n)}).\]
We now describe this $S_n$-action on $H^*_T(G/B)$ explicitly by giving a formula for the
action of a simple transposition on an equivariant Schubert class.   For an arbitrary permutation $u$,
the action of $u$ on $(p_v)_{v \in S_n}$ is obtained by
factoring $u$ into simple transpositions and then inductively applying the formula.  Since the
$S_n$-action on $H^*_T(G/B)$ is well-defined, the result is independent of the factorization of
$u$---though that is not obvious from the formula for the action of a simple transposition!

Our methods in this section consist entirely of elementary combinatorics.  For examples,
the reader may wish to refer to Figures \ref{examples of perm action} and \ref{egs of KT classes},
which contain the calculations
$s_{2,3}[\Omega_{(12)}]=[\Omega_{(12)}]$ and $s_{1,2}[\Omega_{(12)}]=[\Omega_{(12)}]+(t_2-t_1)[\Omega_e]$.

The reader interested in the situation outside of type $A_n$ should note that 
the following proof will apply immediately to general Lie type, once the generalizations of Lemma
\ref{perm lemma}, \ref{moment graph props}, and \ref{localizing schubert classes} and
Corollary \ref{neighboring polys} are given in Section \ref{general type}.

\begin{proposition} \label{eqvt flag repn}
For each Schubert class $[\Omega_w] \in H^*_T(G/B)$ and simple transposition $s_{i,i+1}$, the action
of $s_{i,i+1}$ on $[\Omega_w]$ is given by
\[s_{i,i+1}[\Omega_w] = \left\{ \begin{array}{ll}
[\Omega_w] & \textup{ if } s_{i,i+1}w > w \textup{ and}\\
{[\Omega_w]}+(t_{i+1}-t_{i})[\Omega_{s_{i,i+1}w}] & \textup{ if } s_{i,i+1}w < w. \end{array} \right.\]
\end{proposition}

\begin{proof}
For convenience, write $s_{i,i+1}[\Omega_w] = (q_v)_{v \in S_n}$.
Schubert classes form a basis for $H^*_T(G/B)$ over $\C[t_1, \ldots, t_n]$.
We expand $s_{i,i+1}[\Omega_w] = \sum c_v^w [\Omega_v]$ in terms of this basis 
 and identify the coefficients $c_v^w$.   To begin,
the degree of $[\Omega_v]$ is $\ell(v)$ 
and the degree of $s_{i,i+1}[\Omega_w]$ is $\ell(w)$ by 
Lemma \ref{localizing schubert classes}.\ref{degree}. 
This means $c_v^w$ is zero  if $\ell(v) > \ell(w)$.

{\bf{ Case 1:}} $s_{i,i+1}w > w$.  We show that if $\ell(v) \leq \ell(w)$
then either $q_v = 0$ or $v=w$. 
Suppose $v$ satisfies $q_v \neq 0$ and $\ell(v) \leq \ell(w)$.
Since $q_v \neq 0$ we know $s_{i,i+1}v \geq w$. 
Recall that $\ell(s_{i,i+1}v) = \ell(v) \pm 1$.  If $\ell(s_{i,i+1}v) = \ell(w)$ then in fact $s_{i,i+1}v = w$, which contradicts $s_{i,i+1}w>w$.  Consequently $\ell(s_{i,i+1}v)=\ell(w)+1=\ell(v)+1$.  
Lemma \ref{localizing schubert classes}.\ref{degree} implies
that $s_{i,i+1}v \rightarrow w$ 
is an edge in the moment graph, say labeled $t_j-t_k$.
If  $s_{jk} \neq s_{i,i+1}$  
then $s_{i,i+1}^2v > s_{i,i+1}v$ by Lemma
\ref{perm lemma}.\ref{transposing reflection}.  This contradicts $s_{i,i+1}v>v$.
We conclude that $s_{jk} = s_{i,i+1}$ and so $v=w$.
By definition
$q_w=s_{i,i+1}p_{s_{i,i+1}w}^w$.
By Corollary \ref{neighboring polys}, we know $s_{i,i+1}p_{s_{i,i+1}w}^w=p_w^w$.
This means $c_w^w=1$ and $c_v^w=0$ for all other $v$, 
namely $s_{i,i+1}[\Omega_w]=[\Omega_w]$.

{\bf Case 2:} $s_{i,i+1}w<w$. First we show that $[\Omega_w]+(t_{i+1}-t_i)[\Omega_{s_{i,i+1}w}]$ 
agrees with $s_{i,i+1}[\Omega_w]$ in entries $w$ and $s_{i,i+1}w$.  
We know that $q_{w}=s_{i,i+1}p_{s_{i,i+1}w}^w$ and $q_{s_{i,i+1}w}=s_{i,i+1}p^w_w$.  
 The Knutson-Tao construction shows that $p_{s_{i,i+1}w}^w=0$.
  Corollary \ref{neighboring polys} says $s_{i,i+1}p_w^w = 
(t_{i+1}-t_{i})p_{s_{i,i+1}w}^{s_{i,i+1}w}$ and   
$p_w^w = -(t_{i+1}-t_i)p_{w}^{s_{i,i+1}w}$. 
It follows that $p_u^w+(t_{i+1}-t_i)p_u^{s_{i,i+1}w}=q_u$ when $u=w$ or
$u=s_{i,i+1}w$.

Suppose $v \neq s_{i,i+1}w$ is a fixed point with $q_v \neq 0$
 and $\ell(v) \leq \ell(w)$.  If $q_v \neq 0$ then $s_{i,i+1}v > w$.  The length restrictions $\ell(v) \leq \ell(w)$ and $\ell(s_{i,i+1}v)> \ell(w)$ together with $\ell(s_{i,i+1}v) = \ell(v) \pm 1$ imply that $\ell(v) = \ell(w)$.  This means the four fixed points $v,w,s_{i,i+1}v,s_{i,i+1}w$ form the fragment of the moment graph shown in Figure \ref{moment graph fragment}.  (The edge between $v$ and $s_{i,i+1}w$ exists and is labeled $t_j-t_k$ since $v = s_{i,i+1}s_{j',k'}s_{i,i+1}w=s_{jk}w$ by virtue of the other three edges.)
\begin{figure}[h]
\begin{picture}(350,34)(0,-19)
\put(175,10){\circle*{5}}
\put(195,0){\circle*{5}}
\put(155,0){\circle*{5}}
\put(175,-10){\circle*{5}}

\put(175,10){\line(-2,-1){20}}
\multiput(175,10)(2,-1){10}{\circle*{1}}
\put(175,-10){\line(-2,1){20}}
\multiput(175,-10)(2,1){10}{\circle*{1}}

\put(180,12){$s_{i,i+1}v$}
\put(200,-2){$v$}
\put(140,-2){$w$}
\put(160,-19){$s_{i,i+1}w$}

\put(90,6){$t_{j'}-t_{k'}$}
\put(130,8){\vector(1,0){35}}

\put(225,-10){$t_j-t_k$}
\put(220,-8){\vector(-1,0){36}}
\end{picture}
\caption{Fragment of moment graph (angles may not be to scale)} \label{moment graph fragment}
\end{figure}
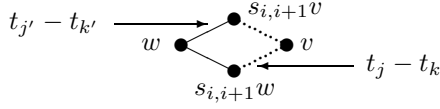
 
We have shown $q_{v} = p_{v}^w+(t_{i+1}-t_i)p_{v}^{s_{i,i+1}w}$ when $v=w$ or $v=s_{i,i+1}w$ and now show it for each $v \neq w$ such that $\ell(v) = \ell(w)$ and $v > s_{i,i+1}w$.
By construction, we have $s_{j'k'}=s_{i,i+1}s_{jk}s_{i,i+1}$
and $s_{i,i+1}(t_{j'}-t_{k'}) = t_j-t_k$.
Lemma \ref{localizing schubert classes}.\ref{neighbors} identifies 
$p_{s_{i,i+1}v}^w$ (respectively $p_{v}^{s_{i,i+1}w}$) as the
product of the monomials in $\Inv{s_{i,i+1}v}$ except $t_{j'}-t_{k'}$
(respectively $\Inv{v}$ and $t_j-t_k$).  
Note that $\Inv{s_{i,i+1}v}$
is exactly $\{t_i-t_{i+1}\} \cup s_{i,i+1}\Inv{v}$
by Lemma \ref{moment graph props}.\ref{up neighbor}.
We conclude $q_{v}=s_{i,i+1}p_{s_{i,i+1}v}^w = (t_{i+1}-t_i)
p_{v}^{s_{i,i+1}w}$.  
Since $p_{v}^w=0$ unless $s_{jk}=s_{i,i+1}$, we have $q_{v}= 
p_{v}^w + (t_{i+1}-t_i)p_{v}^{s_{i,i+1}w}$.
\end{proof}

\subsection{The $S_n$-action on $H^*_T(\Cy{w})$ and $H^*(\Cy{w})$} \label{decomposing repns}

The group  $S_n$ acts on $H^*_T(\Cy{w})$ and on $H^*(\Cy{w})$ by restricting the 
$S_n$-action on $H^*_T(G/B)$.  This section contains the main results of this paper: 
the representation $H^*_T(\Cy{w})$ is not trivial over $\C$ but is trivial
over $\C[t_1, \ldots, t_n]$, and the representation $H^*(\Cy{w})$ is trivial.  (These results hold
for $G/B$ since $G/B = \Cy{w_0}$ when $w_0$ is the longest permutation.)

\begin{lemma} \label{action of longer perms}
For each $s \in S_n$, the action of $s$ on the Schubert class 
$[\Omega_w] \in H^*_T(G/B)$ is given by
\[s \cdot [\Omega_w] = [\Omega_w] + \sum_{v < w} c_v [\Omega_v] \]
for some polynomials $c_v \in \C[t_1,\ldots,t_n]$ with $\deg(c_v) = \ell(w)-\ell(v)$.
\end{lemma}

\begin{proof}
If $s=s_{i,i+1}$ then the claim holds by inspection of the formula in Proposition 
\ref{eqvt flag repn}.
The permutation $s$ can be factored as a product of simple transpositions $s_{i,i+1}$.  
Since $s_{i,i+1}$ preserves degree, induction on the length of $s$ completes the proof.
\end{proof}

\begin{proposition} \label{general simple transposition formula}
The group $S_n$ acts on $H^*_T(\Cy{w})$ by the rule that for each 
$s_{i,i+1} \in S_n$ and Knutson-Tao class $[\Omega_v]_{\Cy{w}} \in H^*_T(\Cy{w})$, we have
\[s_{i,i+1}[\Omega_v]_{\Cy{w}} = \left\{ \begin{array}{ll}
[\Omega_v]_{\Cy{w}} & \textup{ if } s_{i,i+1}v > v  \textup{ and}\\
{[\Omega_v]_{\Cy{w}}}+(t_{i+1}-t_i)[\Omega_{s_{i,i+1}v}]_{\Cy{w}} & \textup{ if } s_{i,i+1}v < v. \end{array} \right.\]
\end{proposition}

\begin{proof}
The restriction map $\iota: H^*_T(G/B) \rightarrow H^*_T(\Cy{w})$ given by 
$\iota \left( (p_u)_{u \in S_n} \right) = (p_u)_{[u] \in [S_n] \cap \Cy{w}}$ is a $\C[t_1,\ldots,t_n]$-module 
homomorphism.  For each $s \in S_n$, we know that $s[\Omega_v] = [\Omega_v]+
\sum_{u < v} c_u [\Omega_u]$ by Lemma \ref{action of longer perms}.
We obtain a well-defined $S_n$-action on $H^*_T(\Cy{w})$
by the rule that $s \cdot [\Omega_v]_{\Cy{w}}=\iota(s[\Omega_v]))$.  This action
 restricts the $S_n$-representation on $H^*_T(G/B)$ in Proposition \ref{eqvt flag repn}.
\end{proof}

Note that this representation is not trivial as a $\C[S_n]$-module.  For instance, 
the submodule $\C[t_1, \ldots, t_n] 
[\Omega_e]_{\Cy{w}}$ of $H^*_T(\Cy{w})$ is invariant under $S_n$.  In fact,
it is isomorphic to the $S_n$-algebra
induced by the standard $S_n$-action on $\C t_1 \oplus \ldots \oplus \C t_n$, which is not
trivial.  The next theorem shows this is essentially the only way in which these
representations are not trivial.  Write $1^d$ for the trivial representation in degree $d$.

\begin{theorem}
The $S_n$-representation $H^*_T(\Cy{w})$ is isomorphic to $\bigoplus_{[v] \in [S_n] \cap \Cy{w}}
1^{\ell(v)}$ as a graded twisted $\C[t_1, \ldots, t_n]$-module.
\end{theorem}

\begin{proof}
The $S_n$-action preserves the degree of each class $[\Omega_v]_{\Cy{w}} \in H^*_T(\Cy{w})$
as well as the degree of each polynomial $p \in \C[t_1, \ldots, t_n]$.  If $s \in S_n$ then we can write $s[\Omega_v]_{\Cy{w}} = [\Omega_v]_{\Cy{w}} + \sum_{u < v} c_u [\Omega_u]_{\Cy{w}}$
for some $c_u \in \C[t_1,\ldots,t_n]$ by Lemma \ref{action of longer perms}.  Let  
\[ [Y_v]_{\Cy{w}} = \displaystyle \frac{\sum_{s \in S_n} s[\Omega_v]_{\Cy{w}}}{|S_n|}\] 
so in particular $[Y_v]_{\Cy{w}}$ is $S_n$-invariant and has degree $\ell(v)$.  

We show that the $[Y_v]_{\Cy{w}}$ generate $H^*_T(\Cy{w})$.  There are polynomials $q_u$ with
$[Y_v]_{\Cy{w}} = [\Omega_v]_{\Cy{w}}+\sum_{u < v} q_u [\Omega_u]_{\Cy{w}}$ by construction.
It follows that $[Y_e]_{\Cy{w}} = [\Omega_e]_{\Cy{w}}$.  This proves the inductive hypothesis that if $v$ has length at most $k$ then the Schubert class $[X_v]_{\Cy{w}}$ can be written in terms of the $[Y_u]_{\Cy{w}}$ for $u \leq v$.  If $v$ has length $k+1$ then $[\Omega_v]_{\Cy{w}} = [Y_v]_{\Cy{w}} - \sum_{u < v} q_u [\Omega_u]_{\Cy{w}}$ which can by induction be written in terms of the $[Y_u]_{\Cy{w}}$ for $u \leq v$.  It follows that the Schubert classes, a known basis for $H^*_T(\Cy{w})$, are generated by the $[Y_v]_{\Cy{w}}$.  So $H^*_T(\Cy{w})$ decomposes as $|\Cy{w} \cap [S_n]|$ distinct $S_n$-modules,
each eigenspaces over $\C[t_1, \ldots, t_n][S_n]$.
\end{proof}

\begin{corollary}
Let $\C[t_1, \ldots, t_n]$ denote the $S_n$-algebra induced by the standard
$S_n$-representation on $\C t_1 \oplus \cdots \oplus \C t_n$.  Then
$H^*_T(\Cy{w})$ is isomorphic to $\bigoplus_{[v] \in [S_n] \cap \Cy{w}}
1^{\ell(v)} \otimes \C[t_1, \ldots, t_n]$ as a graded $\C[S_n]$-module.  
\end{corollary}

The $S_n$-action on $H^*_T(\Cy{w})$ gives rise to a $S_n$-action 
on $H^*(\Cy{w})$.

\begin{theorem}
Consider $H^*(\Cy{w})$ as a $\C[S_n]$-module with the $S_n$-action inherited from $H^*_T(\Cy{w})$.
Then $H^*(\Cy{w})$ is isomorphic to $\bigoplus_{[v] \in [S_n] \cap \Cy{w}} 1^{\ell(v)}$ as a graded 
$\C[S_n]$-module.
\end{theorem}

\begin{proof}
The ring $H^*(\Cy{w})$ is isomorphic to $H^*_T(\Cy{w})/ \langle t_1, \ldots, t_n \rangle$ (see
\cite[Equation 1.2.4]{GKM}).
Consequently, the $S_n$-action on $H^*_T(\Cy{w})$ induces an $S_n$-action
on $H^*(\Cy{w})$ given by $s_{i,i+1}[\Omega_v]_{\Cy{w}} = [\Omega_v]_{\Cy{w}}$ for each
 $s_{i,i+1} \in S_n$ and
$[v] \in \Cy{w} \cap [S_n]$.  This representation is trivial.
\end{proof}

\subsection{Divided difference operators}

This action of $S_n$ on $H^*_T(\Cy{w})$ gives rise to divided difference operators on the
equivariant cohomology ring.

\begin{definition}
For each $i=1,\ldots,n-1$, the $i^{th}$ (left) divided difference operator $D_i$ is defined by
\[D_i (p) = \frac{p - s_{i,i+1} \cdot p}{t_i-t_{i+1}}\]
for each $p \in H^*_T(\Cy{w})$.
\end{definition}

Note that the left divided difference operator uses the $S_n$-action on equivariant
classes rather than the $S_n$-action on the polynomials obtained by localizing
equivariant classes.  

We prove that the divided difference operators are well-defined by using the following 
stronger result.  Let $\partial_i: \C[t_1, \ldots,t_n] \rightarrow \C[t_1,\ldots,t_n]$ denote the
(ordinary) divided difference operator defined by $\partial_i(p) = \frac{p-s_ip}{t_i-t_{i+1}}$.

\begin{proposition}
Let $p \in H^*_T(\Cy{w})$ be an equivariant class whose expansion in terms of the
basis of Knutson-Tao classes is $p = \sum c_v [\Omega_v]$.  Then
\[D_i(p) = \sum \partial_i(c_v)[\Omega_v] + \sum_{v: s_iv<v} s_i(c_v) [\Omega_{s_iv}].\]
\end{proposition}

\begin{proof}
The $S_n$-action on $H^*_T(\Cy{w})$ is $\C$-linear, so 
$D_i(p) = \sum D_i(c_v[\Omega_v])$.  Proposition \ref{general simple transposition formula} 
shows that $D_i(c_v[\Omega_v])$ is $\partial_i(c_v)[\Omega_v]$ if $s_i v > v$ and is 
$\partial_i(c_v)[\Omega_v]+s_i(c_v)[\Omega_{s_iv}]$ if $s_iv<v$. 
\end{proof}

In particular $D_i(p)$ is an element of $H^*_T(\Cy{w})$ whenever $p$ is in $H^*_T(\Cy{w})$.

Every Schubert class of $G/B$ can be obtained by performing a sequence of divided difference operators on the class of the longest permutation $[\Cy{w_0}]$.  This is not true for general Schubert varieties.  For instance, no sequence of divided difference operators performed on the highest class $[\Cy{s_1s_2}]$ of the Schubert variety $X_{s_1s_2}$ gives the class $[\Cy{s_2}]$.

The divided difference operator $D_i$ is not the same as the operator defined by Bernstein-Gelfand-Gelfand and Demazure \cite{BGG}, \cite{D}, though the formulas are similar.  Another divided difference operator was defined by Kostant-Kumar in \cite{KK1}: for each $p \in H^*_T(G/B)$ and each $v \in S_n$ let
\[(\partial_i p)(v) = \frac{p(v) - p(vs_i)}{-v(t_{i} - t_{i+1})}.\]
Arabia proved in \cite{A} that the operator $\partial_i$ is the same morphism as the Bernstein-Gelfand-Gelfand/Demazure divided difference operator.  The formulas differ because Bernstein-Gelfand-Gelfand/Demazure use a different presentation of the ring $H^*_T(G/B)$ than Kostant-Kumar (and we) do.

\section{General Lie type} \label{general type}

We now describe these results in arbitrary Lie type.  Our exposition is brief: we assume our
reader is familiar with the general theory and only indicate the proofs whose generalization
is not immediate.

In this section, $G$ is a complex reductive linear algebraic group,
$B$ a Borel subgroup, and $T$ a maximal torus contained in $B$.  The full flag variety is $G/B$.
We write an element of the flag variety as $[g]$.  
The Weyl group is  $W = N(T)/T$.  
The set of simple reflections in $W$ is written $S$ and the set
of all reflections is $R = \bigcup_{w \in W} wSw^{-1}$.
The roots are denoted $\Phi$, with positive roots $\Phi^+$ and negative roots $\Phi^-$.  The
elements in $R$ are the reflections associated to the positive roots.  We write $s_{\alpha} \in R$
to denote the reflection associated to $\alpha \in \Phi^+$.

There is a partial order on $W$ called the Bruhat order and defined by the condition that
$w \geq v$ if and only if $\overline{[Bw]} \supseteq [Bv]$.  This is equivalent to the condition that
there is a reduced factorization $w = s_{i_1} \cdots s_{i_k}$ with each $s_{i_j} \in S$ such that
$v$ is the (ordered) 
product of a substring of the $s_{i_j}$.  The length of $w$ is the minimal $k$ required
to factor $w = s_{i_1} \cdots s_{i_k}$ as a product of simple transpositions $s_{i_j} \in S$.  
The length of $w$ is denoted $\ell(w)=k$.

The next few results are combinatorial properties of $W$.  This next proposition complements
\cite[Lemma 2.4]{BGG}.

\begin{proposition}
Suppose $s_{\alpha} \in R$ and $w \in W$ satisfy $\ell(s_{\alpha}w) = \ell(w)+1$.
Then
\begin{enumerate}
\item $s_{\alpha}w > w$; and
\item if $s_i$ is a simple reflection with $s_i \neq s_{\alpha}$ and $s_iw > w$,
 then $s_is_{\alpha}w > s_{\alpha}w$.
\end{enumerate}
\end{proposition}

\begin{proof}
The first part follows because either $s_{\alpha}w > w$ or $s_{\alpha}w<w$, and the Bruhat
order respects length.  For the second
part, suppose $s_i s_{\alpha}w < s_{\alpha}w$.  By the exchange property \cite[IV.1.5]{Bou} or \cite[Corollary 1.4.4]{BB}, there is a 
reduced expression for $s_{\alpha}w$ that begins with $s_i$, say $s_i s_{i_1} \cdots s_{i_k}$.
Since $s_{\alpha}w > w$ and $\ell(s_{\alpha}w) = \ell(w)+1$, we can obtain a reduced 
factorization of $w$ by erasing one of the simple transpositions in the string
$s_i s_{i_1} \cdots s_{i_k}$.  If it were not $s_i$ then
$s_iw < w$.  So $s_i s_{\alpha}w > s_{\alpha}w$ unless $s_{\alpha} = s_i$.
\end{proof}

The following set is well-known in the literature (e.g. \cite[VI.1.6, Proposition 17]{Bou}) but seems not to have a concise name.

\begin{definition}
Given $w \in W$, the inversions corresponding to $w$ are the roots
\[\Inv{w} = \{\alpha \in \Phi^+: w^{-1}\alpha \in \Phi^-\} = \Phi^+ \cap w \Phi^-.\]
\end{definition}

For instance, each simple transposition $s_i$ has a unique inversion $\Inv{s_i} = \{\alpha_i\}$.
There are $\ell(w)$ elements in the set $\Inv{w}$.

\begin{proposition} \label{simple edge}
If $s_iw > w$ then $\Inv{s_iw} = s_i \Inv{w} \cup \{\alpha_i\}$.
\end{proposition}

\begin{proof}
The definition of $\Inv{w}$ implies that $(w^{-1}s_i) s_i \Inv{w} \subseteq \Phi^-$.  
By the exchange property \cite[IV.1.5]{Bou} or \cite[Corollary 1.4.4]{BB}, there is a reduced factorization for $s_iw$ beginning
with $s_i$.  It is a small exercise (proven in, e.g., \cite[Corollary 10.2.C]{H}) to see that this implies 
that $(s_iw)^{-1} (\alpha_i)$ is negative.  We conclude that $\alpha_i \in \Inv{s_iw}$ and that
$\alpha_i \not \in \Inv{w}$.
 The simple reflection $s_i$ negates $\alpha_i$
and permutes the other positive roots, so $s_i \Inv{w} \subseteq \Phi^+$.  We have shown
that the set
$s_i \Inv{w} \cup \{\alpha_i\}$ consists of distinct roots and is contained in $\Inv{s_iw}$.  Since both
sets contain $\ell(w)+1$ roots, they are in fact equal.  
\end{proof}

\begin{proposition}
If $s_{\alpha}w > w$ and $\ell(s_{\alpha}w) = \ell(w)+1$ then
\[\Inv{s_{\alpha}w} = \{\alpha\} \cup \left( \Phi^+ \cap s_{\alpha} \Inv{w} \right)  \cup  
 \left( \Inv{w} \cap s_{\alpha} \Phi^- \right).\]
\end{proposition}

\begin{proof}
The root $\alpha$ must be in exactly one of $\Inv{w}$ and $\Inv{s_{\alpha}w}$.  Since
$s_{\alpha}w>w$, we have $\alpha \in \Inv{s_{\alpha}w}$ (see \cite[Theorem F4]{C} or
Proposition \ref{general moment graph}.\ref{edges in general}.).
Note that if $\beta \in \Inv{w}$ then $(s_{\alpha}w)^{-1} s_{\alpha}(\beta)$ is a negative root.
If $s_{\alpha}(\beta) \in \Phi^+$ then $s_{\alpha}(\beta) \in \Inv{s_{\alpha}w}$.  Suppose
instead that $s_{\alpha}(\beta) \in \Phi^-$.  There is a positive
integer $c_{\beta}$ such that $s_{\alpha}(\beta) = \beta - c_{\beta}\alpha$ by definition
of reflections.  On the one hand, we know that $\beta=s_{\alpha}(\beta) + c_{\beta}\alpha$
 is a positive root.  On the other hand,
we know that $w^{-1}s_{\alpha} (s_{\alpha}(\beta) + c_{\beta}\alpha) = 
w^{-1}(\beta) + c_{\beta} w^{-1}(-\alpha)$.  We conclude that $\beta$ is in $\Inv{s_{\alpha}w}$,
since $w^{-1}(\beta)$ and $w^{-1}(-\alpha)$ are both
negative and $c_{\beta}$ is positive.  

We have shown that the map sending $\beta \in \Inv{w}$ to
$\left\{ \begin{array}{ll} s_{\alpha}(\beta) & \textup{ if }s_{\alpha}(\beta) \in \Phi^+ \textup{ and}
\\ \beta & \textup{ else} \end{array} \right.$ has image in $\Inv{s_{\alpha}w}$.  It is an injection
because it is invertible.  By comparing cardinalities, we know that $\Inv{s_{\alpha}w}$ consists
of the image of this map together with exactly one other root.  This proves the claim.
\end{proof}

\begin{corollary} \label{general s_alpha w}
If $s_{\alpha}w > w$ and $\ell(s_{\alpha}w) = \ell(w)+1$ then
\[\Inv{s_{\alpha}w} = \{\alpha\} \cup s_{\alpha} \Inv{w} \mod(\alpha) \textup{  (with multiplicity).}\]
\end{corollary}

\begin{proof}
If $\beta$ in $\Inv{w}$ and $s_{\alpha}(\beta) \not \in \Inv{s_{\alpha}w}$ then by the previous
proposition, there is a positive integer $c_{\beta}$ with $\beta = s_{\alpha}(\beta)
+ c_{\beta} \alpha \in \Inv{s_{\alpha}w}$.  In this case $\beta \equiv s_{\alpha}(\beta) \mod \alpha$.
\end{proof}

J.~ Carrell described the moment graph for Schubert varieties in 
$G/B$ in \cite[Theorem F]{C}.  We give his result
here for the convenience of the reader.  The root subgroup corresponding to $\alpha$
is written $U_{\alpha}$.

\begin{proposition} \label{general moment graph}
\begin{enumerate}
\item The vertices for the moment graph for $G/B$ are exactly the flags $[w]$ for $w \in W$.
\item \label{edges in general}
	There is an edge between $[v]$ and $[w]$ in the moment graph if and only if
	$v^{-1}w \in R$.  If $v^{-1}w = s_{\alpha}$ then the edge is labeled $\alpha$ and
	corresponds to the one-dimensional $T$-orbit $[U_{\alpha}w] \cup [U_{\alpha}v]$.
	The edge between $[s_{\alpha}v]$ and $[v]$ 
	is directed from the element of greater length to the element
	of lesser length.
\item The edges directed out of $[w]$ are labeled exactly by the roots in $\Inv{w}$.
\item The moment graph for the Schubert variety $\Cy{w} = \overline{[Bw]}$ is the
	subgraph induced by the elements $[v] \in [W] \cap \Cy{w}$ from the moment
	graph for $G/B$.
\end{enumerate}
\end{proposition}

Write $H^*_T(\textup{pt})$ as the polynomial ring in the simple roots $\C[\alpha_1, \ldots,
\alpha_n]$.  The Weyl group acts on $H^*_T(\textup{pt})$ by extending the coadjoint
$W$-action on $\Phi$, namely
 $w \in W$ acts on $p(\alpha_1, \ldots, \alpha_n) \in H^*_T(\textup{pt})$ by
$w \cdot p(\alpha_1, \ldots, \alpha_n) = p(w(\alpha_1), \ldots, w(\alpha_n))$.  
The equivariant cohomology of each Schubert variety $H^*_T(\Cy{w})$ is a 
module over $H^*_T(\textup{pt})$.

In particular, let $[\Omega_v]$ be the class in $H^*_T(G/B)$ determined by the closure
$\overline{[B^-vB]}$ in $G/B$.  The previous results imply that $[\Omega_v]$ is
the Knutson-Tao class for $v$.

\begin{proposition}  Let $v \in W$ and write $[\Omega_{v}] \in H^*_T(G/B)$
as $[\Omega_{v}]=(p^v_u)_{u \in W}$.
\begin{enumerate}
\item If $u \not \geq v$ then $p^v_u = 0$.  If $u \geq v$ then $p^v_u$ has degree $\ell(v)$.
\item The localization of $[\Omega_v]$ to $v$ is $p^v_v = \prod_{\alpha \in \Inv{v}} \alpha$.
\item The class $[\Omega_{v}]$ is the unique Knutson-Tao class for $v$ in $H^*_T(G/B)$.
\item If $u > v$ satisfies $u = s_{\alpha}v$ and $\ell(u)=\ell(v)+1$ then
\[p^v_u = \frac{1}{\alpha} \prod_{\beta \in \Inv{u}} \beta.\]
\item If $s_i \in S$ satisfies $s_iv > v$ then
\[s_i p^v_{s_iv} = p^v_v = \frac{s_ip^{s_iv}_{s_iv}}{-\alpha_i}.\]
\end{enumerate}
\end{proposition}

\begin{proof}
As in type $A_n$, the closure $\overline{[B^-vB]}$
contains exactly the $T$-fixed points $[w]$ with $w \geq v$.  The normal space
to $\overline{[B^-vB]}$ at $v$ is given by $[Bv]$, whose torus weights are indexed by
$\Inv{v}$.  So, the 
proofs of Parts 1, 2, and 3 are identical to those in Lemma \ref{localizing schubert classes}.   
The proof of the other parts is closely parallel to those for type $A_n$.
Write $q =
\frac{1}{\alpha} \prod_{\beta \in \Inv{u}} \beta$.  The polynomial $p^v_u \in \langle q \rangle$ 
by the GKM conditions, so $p^v_u = cq$ for some (complex) constant $c$.  
Corollary \ref{general s_alpha w} shows $q - p^v_v \in \langle \alpha \rangle$, so $c = 1$.
This proves Part 4.  Part 5 follows directly from Part 2 and Proposition \ref{simple edge}.
\end{proof}

We now describe a group action of $W$ on $H^*_T(G/B)$ induced by the action of $u \in W$ on $[v] \in [W]$ given by $u \cdot [v] = [u^{-1}v]$.  Unlike the case when $G = GL_n(\C)$, this group action is not well-defined on the entire flag variety, since $u \in W = N(T)/T$ is defined only up to $T$.  However, the action does induce a graph automorphism of the moment graph of $G/B$.  This action arises from the action of $u \in W$ on $ET \times^T G/B$ given by $u \cdot (e,[g]) = (eu, [u^{-1}g])$ just as in Corollary \ref{W-action}.

\begin{proposition}
The group $W$ acts on $H^*_T(G/B)$ according to the rule that if $u \in W$ and 
$p = (p_v)_{v \in W} \in H^*_T(G/B)$ then for each $v \in W$
\[ \left( u \cdot p \right)_{v \in W} = u \cdot p_{u^{-1}v}\left(\alpha_1,
\ldots, \alpha_n \right) = p_{u^{-1}v}\left(u(\alpha_1),
\ldots, u(\alpha_n) \right).\]
\end{proposition}

\begin{proof}
We give the main step of the proof in detail; the rest 
is exactly as in Corollary \ref{W-action}. The action of 
$u^{-1} \in W$ on $G/B$ given by $u^{-1} \cdot [g] = [u^{-1}g]$
is defined on each one-orbit $[U_{\alpha}v]$ since each $T$-one-orbit is $T$-closed.
The action satisfies $u^{-1} \cdot [U_{\alpha}v] = [U_{u^{-1}(\alpha)}u^{-1}v]$.  
If $u^{-1}(\alpha) \in \Phi^+$ then this is 
a one-orbit between $[u^{-1}v]$ and $[s_{u^{-1}(\alpha)}u^{-1}v]$ in $G/B$.  
If $u^{-1}(\alpha) \in \Phi^-$ then
$[U_{-u^{-1}(\alpha)}s_{u^{-1}(\alpha)}u^{-1}v]$ is a one-orbit between 
$[s_{u^{-1}(\alpha)}u^{-1}v]$ and $[u^{-1}v]$
in $G/B$.   In either case, there is an edge between $v$ and $s_{\alpha}v$ if and only if
there is an edge between $u^{-1}v$ and $u^{-1}s_{\alpha} v = s_{u^{-1}(\alpha)} u^{-1}v$.
\end{proof}

The propositions included in this section are exactly those needed in Proposition
\ref{eqvt flag repn}.  
The proof of Proposition \ref{eqvt flag repn} applies in our more general setting once the 
appropriate
notational changes are made: $\alpha_i$ for $t_i - t_{i+1}$, $s_i$ for $s_{i,i+1}$, 
$\alpha$ for $t_j-t_k$, and $s_{\alpha}$ for $s_{jk}$.  We give the statement of the general 
theorem here.

\begin{theorem}
For each Schubert class $[\Omega_w] \in H^*_T(G/B)$ and simple transposition $s_i$, the 
action of $s_i$ on $[\Omega_w]$ is given by 
\[s_{i}[\Omega_w] = \left\{ \begin{array}{ll}
[\Omega_w] & \textup{ if } s_{i}w > w \textup{ and}\\
{[\Omega_w]}- \alpha_i[\Omega_{s_{i}w}] & \textup{ if } s_{i}w < w. \end{array} \right.\]
\end{theorem}

The restriction of this action to $H^*_T(\Cy{w})$ and the analysis of the $W$-action on
$H^*_T(\Cy{w})$ and $H^*(\Cy{w})$ is also identical to the case of type $A_n$.  The statement
is:

\begin{theorem} Let $\Cy{w}$ be a Schubert variety in $G/B$.  The trivial representation
of $W$ in degree $d$ is denoted $1^d$, 
and the $W$-algebra induced by the coadjoint action on $\Phi$
is denoted $\C[\alpha_1, \ldots, \alpha_n]$.
\begin{enumerate}
\item $H^*_T(\Cy{w})$ is isomorphic to $\bigoplus_{[v] \in [W] \cap \Cy{w}}
1^{\ell(v)} \otimes \C[\alpha_1, \ldots, \alpha_n]$ as a graded twisted $\C[W]$-module;
\item $H^*_T(\Cy{w})$ is isomorphic to $\bigoplus_{[v] \in [W] \cap \Cy{w}}
1^{\ell(v)}$ as a graded twisted module over $\C[\alpha_1,\ldots,\alpha_n][W]$;
\item and $H^*(\Cy{w})$ is isomorphic to $\bigoplus_{[v] \in [W] \cap \Cy{w}}
1^{\ell(v)}$ as a graded $\C[W]$-module.
\end{enumerate}
\end{theorem}

The left $W$-action also gives rise to (left) divided difference operators on the equivariant
cohomology of Schubert varieties in general type.

\begin{definition}
For each $s_i \in W$, the divided difference operator $D_i$ is defined by
\[D_i(p) = \frac{p - s_i \cdot p}{\alpha_i}\]
for each $p \in H^*_T(\Cy{w})$.
\end{definition}

As before, there is an explicit formula for the action of $D_i$.  Let $\partial_i:
\C[\alpha_1,\ldots,\alpha_n] \rightarrow \C[\alpha_1,\ldots,\alpha_n]$ denote the
$i^{th}$ (ordinary) divided difference operator.

\begin{proposition}
If $p = \sum c_v [\Omega_v]$ is an equivariant class in $H^*_T(\Cy{w})$ then
\[D_i(p) =  \sum \partial_i(c_v)[\Omega_v] + \sum_{v: s_iv<v} s_i(c_v) [\Omega_{s_iv}].\]
\end{proposition}

\end{document}